\documentstyle[10pt]{article}
\input amssymb.sty
\input epsf

\newtheorem{theorem}{Theorem}[section]
\newtheorem{lemma}[theorem]{Lemma}
\newtheorem{proposition}[theorem]{Proposition}

\newtheorem{conjecture}[theorem]{Conjecture}

\newtheorem{question}[theorem]{Question}

\begin{document}
\newcommand{\Z}{\Bbb Z}
\newcommand{\R}{\Bbb R}
\newcommand{\Q}{\Bbb Q}
\newcommand{{\C}}{\Bbb C}
\newcommand{\lra}{\longrightarrow}
\newcommand{\lms}{\longmapsto}
 
\begin{titlepage}
\title{ Geometry of the trilogarithm and the motivic Lie algebra of a field}
\author{A.B. Goncharov }

\date{}
\end{titlepage}
\maketitle
\tableofcontents

$$
$$
 
{\bf Summary}. We express explicitly the Aomoto trilogarithm by
classical trilogarithms and
 investigate the algebraic-geometric structures staying behind: 
 different realizations of the weight three motivic
complexes. Applying these results we describe the 
motivic structure  of the Grassmannian tetralogarithm function and determine
 the structure of the motivic Lie coalgebra in  degrees $\leq
4$. Using this we 
give an explicit construction of the Borel regulator map 
$$
r_4: K_7({\C}) \lra \R
$$
which together with the Borel theorem  leads to results about $\zeta_F(4)$. 
 
\section{Introduction}

The classical $n$-logarithm is defined by induction as an integral
$$
Li_n(z):= \int_0^z Li_{n-1}(t)d\log t, \quad Li_1(z)= - \log (1-z)
$$
So it can be written as an $n$-dimensional integral
$$
Li_n(z) = \int_{ 0  \leq 1 - t_1 \leq t_2 \leq ... \leq t_n \leq z}\frac{d t_1}{t_1}\wedge ... \wedge\frac{d t_n}{t_n} 
$$
Aomoto considered [A] more general integrals  where  the differential form $\frac{d t_1}{t_1}\wedge ... \wedge\frac{d t_n}{t_n}$ is integrated over an arbitrary $n$-dimensional real simplex in ${\C}^n$. Let me recall this construction in a more formal setting.

Let $F$ be a field. An
 $n$-simplex in ${\Bbb P}^n(F)$ 
 is  a collection of hyperplanes $L = (L_0,...,L_n)$. It is
nondegenerate if   the intersection of the hyperplanes $L_i$ is empty. 
 A face of a simplex is any nonempty intersection of
hyperplanes from $L$. A pair of simplices is {\it admissible} if $L$ and $M$
have no common faces of the same dimension.
 Now let $F = {\Bbb C}$. Then there is a canonical $n$-form $\omega_L$ in $\Bbb
CP^n$ with logarithmic poles on the hyperplanes $L_i$. If $z_i = 0$ are
homogeneous equations of $L_i$ then 
$$
\omega_L = d\log(z_1/z_0) \wedge
...  \wedge d\log(z_n/z_0)
$$
 Let $\Delta_M$ be an $n$-cycle representing
a generator of $H_n({\Bbb C}{\Bbb P}^n,  M )$. 

The Aomoto $n$-logarithm is a
multivalued function on configurations of admissible pairs of simplices
$(L;M)$ in ${\Bbb C}P^n$ defined as follows:
$$
\Lambda_n(L;M) := \int_{\Delta_M} \omega_L
$$

The classical $n$-logarithm corresponds to a very special pair of
simplices in ${\Bbb P}^n$, see fig.1, so it is a very special case of the Aomoto $n$-logarithm.

\begin{center}
\hspace{4.0cm}
\epsffile{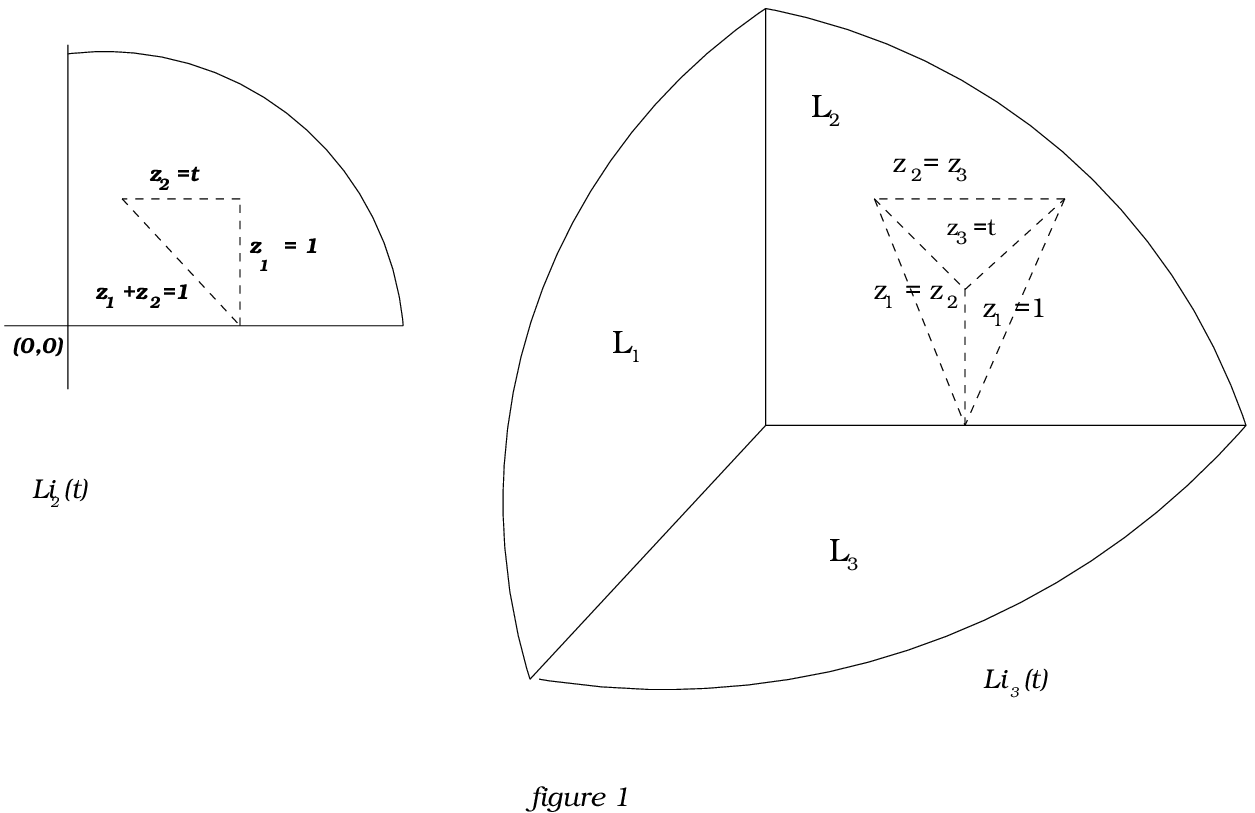}
\end{center}

In the present paper, which  is a continuation of [BGSV1-2] and [G0-3],
we express explicitly the Aomoto trilogarithm by
classical trilogarithms and
 investigate the algebraic-geometric structures staying behind: 
 different realizations of the weight three motivic complexes. 
For the Aomoto dilogarithm a  similar problem was
solved  in [BGSV1-2]. 

The function $\Lambda_3(L;M)$ is defined on configurations, i.e. projective 
equivalence classes  of $4+4$ points in ${\Bbb P}^3$ (vertices of a pair of 
tetrahedra), 
while the classical trilogarithm lives on ${\Bbb P}^1$. 
 To build a bridge between these functions we relate  each of them with the 
Grassmannian trilogarithm ${\cal P}^G_3$ defined on configurations 
of six points 
in ${\Bbb P}^2$ (see section 4.1 for a definition of Grassmannian polylogarithms).
Our main geometric construction, the map $a_3$
defined in s. 3.3, (see its different versions on fig. 3, 8 and 9), 
 permits to go from   configurations 
of $4+4$ points in ${\Bbb P}^3$ to configurations of $6$ points in ${\Bbb P}^2$. 
Then 
we apply the generalized 
cross-ratio map $r_3$ from  [G0-3] to get to ${\Bbb P}^1$. 
The map $a_3$  sheds a new light on the key ansatz 
from  [G0-3] leading to the functional equation for the  classical 
trilogarithm, see fig. 10 and the discussion there.

The  Aomoto $n$-logarithm for $n>3$ {\bf can not} be
expressed by the classical $n$-logarithm. However the 
{\it explicit} relation between the 
Aomoto $n$-logarithm, which is defined on configurations 
of $2(n+1)$ points 
in ${\Bbb P}^{n}$, and the Grassmannian $n$-logarithm, 
which lives on configurations of 
$2n$ 
points in ${\Bbb P}^{n-1}$, should exist for all $n$. In  section 4 we explain 
how such a relation would give an explicit construction of a certain graded 
co-Lie 
algebra $G_{\bullet}(F)$ over ${\Q}$, which should be isomorphic to 
the  Lie coalgebra  of 
the Galois group of the category of mixed Tate motives over a field $F$. The cohomology of this Lie algebra should give the appropriate 
pieces of Quillen's K-theory of the field $F$ modulo torsion. It would be very interesting to relate our approach
 with the work of M. Hanamura and R. MacPherson [HM].

There are several other candidates for the motivic Lie coalgebra, 
see [BK],  [BMS-BGSV], [G5]. However all of them are constructed as  
Hopf algebras, so we get the  Lie coalgebras as the  quotient by the 
decomposable elements.  
Our approach should lead directly to a  Lie coalgebra. The degree $2$ and $3$ 
parts of its standard cochain complex  are precisely the 
Bloch-Suslin complex and the weight three motivic complex defined in [G0-3], 
so 
$G_{\bullet}(F)$ should be the smallest possible realization of the 
motivic Lie algebra.

   In  section $5$ we define the structure of the motivic Lie coalgebra in   degree $4$, i.e. we define a cobracket
$$
G_4(F) \lra G_3(F) \otimes G_1(F) \oplus \Lambda^2 G_2(F)
$$
which satisfies the condition $\delta^2 =0$ in $G_2(F) \otimes \Lambda^2 G_1(F)$.
This together with the previous results of the author provide a description of the Lie subcoalgebra  $G(F)_{\leq 4}$.   

An immediate application of this  is a description of the "fine" (or
motivic) structure  of the Grassmannian tetralogarithm function. In
particular we get an explicit construction of the Borel regulator map 
$$
r_4: K_7({\C}) \lra \R
$$
This together with the famous theorem of Borel leads to results about special values of the Dedekind zeta function at $s=4$ (these results, however, are not sufficient to establish   Zagier's conjecture at $s=4$). 

In the forthcoming paper [G8] we will complete this story by giving an explicit description of the general Beilinson regulator map in  weight $4$. 

Our results partially 
generalize the work [BGSV] from ${\Bbb P}^2$ to ${\Bbb P}^3$. We discuss in s. 4 what remains to be done.  An application  to an explicit construction of the weight four motivic complexes   will be discussed elsewhere.

I am extremely grateful to Herbert Gangl who  wrote a proof of lemma (\ref{3.7'}), helped me to check coefficients in the formulas and   pointed out a lot of errors in a preliminary version of the paper. Finally, I am very much indebted to the  referee who spot a lot of misprints and made many useful remarks.

The results of this paper were obtained in May 1992 during my 
stay in the Max-Planck-Institute (Bonn). The paper was written in MPI later on  
(when I learned how to draw pictures using computer). 
I am very grateful to the MPI for hospitality and support. 

The work was partially supported by the NSF grant DMS-9500010.

\section{The scissors 
congruence groups  of pairs of simplices in ${\Bbb P}^n(F)$}

{\bf 1.The scissors 
congruence groups $A_n(F)$}. 
Let me recall some definitions from [BMS], [BGSV]. It is handier  to work with configurations of points than  with
hyperplanes.  Let us apply the projective duality ${\Bbb P}^n \longrightarrow
\hat {\Bbb P}^n$ which transforms the configuration  of $2(n+1)$ hyperplanes 
$(L_0,...,L_n;M_0,...,M_n)$ to a configuration  of points
$(l_0,...,l_n;m_0,...,m_n)$.   
Abusing notations we will denote it also $(L;M)$, where now $L= (l_0,...,l_n)$ and $M = (m_0,...,m_n)$.

The group $A_n(F)$ is generated by configurations of $2(n+1)$ points
\linebreak $(l_0,...,l_n;m_0,...,m_n)$ in ${\Bbb P}^n(F)$ which are vertices of  admissible
pairs of simplices subject to the following relations:

\vskip 3mm \noindent

1){\it Nondegeneracy}.  $(L;M) = 0$ if $(l_0,...,l_n)$ or $(m_0,...,m_n)$ belong to a
hyperplane.

2) {\it Skew symmetry}. $(\sigma L;M) = (L;\sigma M) = (-1)^{|\sigma|} (L;M)$
for any permutation $\sigma \in S_{n+1}$.

3){\it Additivity}.  For any configuration $(l_0,...,l_{n+1})$ 
$$
\sum_{i=0}^{n+1}(-1)^i (l_0,...,\hat l_i,...,l_{n+1};m_0,...,m_n) =0
$$
if all  the terms are admissible (additivity in $L$).  
A similar
condition is imposed for $(m_0,...,m_{n+1})$ (additivity in $M$).

$\hat 3$){\it Dual additivity}. For any configuration $(l_0,...,l_{n+1})$ 
$$
\sum_{i=0}^{n+1}(-1)^i (l_i\vert l_0,...,\hat l_i,...,l_{n+1};m_0,...,m_n) =0
$$
if all  the terms are admissible,  as well as the similar
condition  is imposed for $(m_0,...,m_{n+1})$. Here 
$(l \vert m_1,...,m_n)$ denotes the  configuration of $n$ points in ${\Bbb P}^{n-1}$ obtained by the projection of the points $m_i$ with the center at the point $l$. 

4) {\it Projective invariance}. $(gL;gM) = (L;M)$ for any $g \in PGL_{n+1}(F)$.

\vskip 3mm \noindent

These  relations reflect  properties of Aomoto polylogarithms. 

The cross-ratio provides a
canonical isomorphism
$$
a_1: A_1(F) \longrightarrow F^{\ast}, \quad a_1:(l_0,l_1;m_0,m_1) \longmapsto 
r(l_0,l_1,m_0,m_1)
$$

{\bf 2. A coproduct on the generic part of $A_n(F)$} ([BMS], [BGSV]). Set $A_0 = {\Bbb Z}$. Let $A^0_n(F) \subset A_n(F)$ be the subgroup generated by pairs of 
simplices in generic position. Let us define  a coproduct $\nu : A^0_n
\longrightarrow \oplus_{k} A_k^0 \otimes A^0_{n-k}$. 
 Set for $k,n-k >0$
$$
 \nu = \oplus \nu_{n-k,k},  \qquad 
\nu_{n-k,k}: A^0_n:  \longrightarrow A_{n-k}^0 \otimes A^0_k,
$$
$$
\nu_{n-k,k}:(l_0,...,l_n;m_0,...,m_n) \longmapsto 
$$
$$
\sum_{I,J}(-1)^{\sigma(I,J)}(l_{i_1}...l_{i_k}\vert l_0,...,\hat
l_{i_1},...,\hat l_{i_k},...,l_n;m_0,m_{j_1},...,m_{j_{n-k}}) 
$$
$$
\otimes (m_{j_1}...m_{j_{n-k}}\vert l_0,
l_{i_1},..., l_{i_k};m_0,...,\hat m_{j_1},...,\hat m_{j_{n-k}},...,m_n) 
$$
where $I:= \{0< i_1 < ... < i_k\}$, $J:= \{0< j_1 < ... < j_{n-k}\}$. Here 
$\sigma(I,J) = sign(I, \bar I) \cdot sign(J, \bar J)$, where 
  $sign(I,\bar I)$ is the 
sign of the permutation $(1,...,n)\rightarrow (I, \bar I)$ ( similarly for 
$sign(J, \bar J)$). Here $\bar I$ is the complement to the set $I$ in $\{0,...,n\}$. 
For example
$$
\nu_{2,1}: (l_0,...,l_3;m_0,...,m_3) \longmapsto 
$$
$$
(l_3\vert l_0,l_1,l_2;m_0,m_2,m_3) \otimes (m_2,m_3\vert l_0,l_3;m_0,m_1)
+ ...
$$

\begin {proposition} 
$(A_{\bullet}(F)^0, \nu)$ is a graded coalgebra.
\end {proposition}

I will need this statement only in the case of 
 degree 3. In this case it is not
hard to deduce it from proposition (\ref{cop}) below.

\begin {proposition} \label{P}
$P: (L;M) \longrightarrow (M;L)$ defines an antiautomorphism of the graded coalgebra
$(A_{\bullet}(F)^0, \nu)$.
\end {proposition}

Proof follows from the definitions.

{\bf 3. Another formula for   $\nu_{1,n-1}$ and $\nu_{n-1,1}$}. 
Unfortunately in the  definition of the coproduct $\nu$ we have to choose
vertices $l_0$   
in $L$ and $m_0$ in $M$ first, so the skew-symmetry is not obvious. In the next
proposition we give another formula for $\nu_{1,n-1}$ and $\nu_{n-1,1}$
which is skew-symmetric from the beginning and  much more convenient.

Let $V_n$ be an $n$-dimensional vector space over a field $F$. Choose a volume form $\omega_n \in det V_n^*$. For any $n$ vectors $l_1,...,l_n$ in $V_n$ set 
$$
\Delta(l_1,...,l_n):= <l_1\wedge ... \wedge l_n, \omega_n >
$$

\begin {proposition} \label{cop}(See figure 2)
$$
\nu_{1,n-1}:(l_0,...,l_n;m_0,...,m_n) \longmapsto 
$$
$$
-  \sum_{i,j =0}^n (-1)^{i+j} \Delta(m_j,l_0,...,\hat
l_i,...,l_n)  \otimes 
(m_j\vert l_0,...,\hat
l_i,...,l_n;m_0,...,\hat m_j,...,m_n)  
$$
$$
\nu_{n-1,1}:(l_0,...,l_n;m_0,...,m_n) \longmapsto 
$$
$$
- \sum_{i,j =0}^n (-1)^{i+j}(l_i\vert l_0,...,\hat 
l_i,...,l_n;m_0,...,\hat m_j,...,m_n)  \otimes
\Delta(l_i,m_0,...,\hat
m_j,...,m_n)   
$$
\end {proposition}

\begin{center}
\hspace{4.0cm}
\epsffile{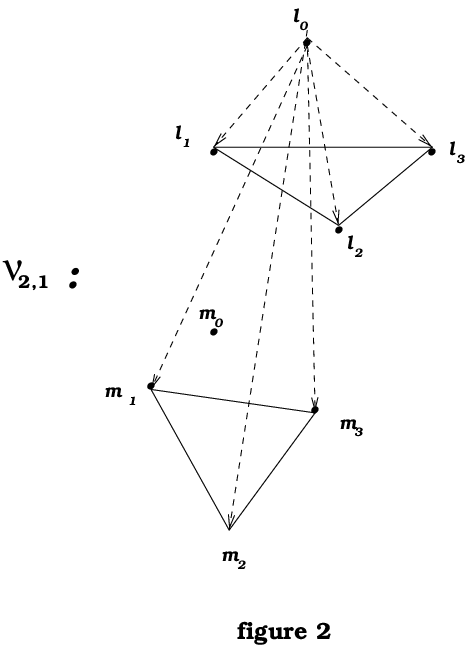}
\end{center}

{\bf Proof}. Let us compute
\begin {equation} \label{a0}
\nu_{n-1,1}(l_0,...,l_n;m_0,...,m_n)  = 
\end {equation}
$$
\sum_{i,j =1}^n (-1)^{i+j}  
(l_i\vert l_0,...,\hat 
l_i,...,l_n;m_0,...,\hat m_j,...,m_n) \otimes r(m_1...\hat
m_j...m_n \vert l_0,l_i;m_0,m_j)
$$
using the  formula  
$$
r(l_1,l_2,l_3,l_4) = \frac{\Delta(l_1,l_3)\Delta(l_2,l_4)}{\Delta(l_1,l_4)\Delta(l_2,l_3)}
$$
 for the cross-ratio. We will get
\begin {equation} \label{a1}
\sum_{i,j =1}^n (-1)^{i+j}  
(l_i\vert l_0,...,\hat 
l_i,...,l_n;m_0,...,\hat m_j,...,m_n) \otimes \Delta(l_0,m_0,...,\hat
m_j,...,m_n)
\end {equation}
\begin {equation} \label{a2}
- \sum_{i,j =1}^n (-1)^{i+j}  
(l_i\vert l_0,...,\hat 
l_i,...,l_n;m_0,...,\hat m_j,...,m_n) \otimes \Delta(l_i,m_0,...,\hat
m_j,...,m_n)
\end {equation}
\begin {equation} \label{a3333}
- \sum_{i,j =1}^n (-1)^{i+j}  
(l_i\vert l_0,...,\hat 
l_i,...,l_n;m_0,...,\hat m_j,...,m_n) \otimes \Delta(l_0,m_1,...,m_n )
\end {equation}
\begin {equation} \label{a4}
+\sum_{i,j =1}^n (-1)^{i+j}  
(l_i\vert l_0,...,\hat 
l_i,...,l_n;m_0,...,\hat m_j,...,m_n) \otimes \Delta(l_i,m_1,...,m_n)
\end {equation}

Applying to (\ref{a1}) dual additivity in $L$ we can rewrite it as 
$$
- \sum_{j =1}^n (-1)^{j}  
(l_0\vert l_1,...,l_n;m_0,...,\hat m_j,...,m_n) \otimes \Delta(l_0,m_0,...\hat
m_j,...,m_n)
$$
(\ref{a2}) is already in  the desired shape. Applying to (\ref{a3333}) first additivity in $M$ and then dual additivity
in $L$ we get
$$
-  
(l_0\vert l_1,...,l_n;m_1,...,m_n) \otimes \Delta(l_0,m_1,...,m_n )
$$
Finally applying additivity in $M$ to (\ref{a4}) we get 
$$
-\sum_{i=1}^n (-1)^{i}  
(l_i\vert l_0,...,\hat 
l_i,...,l_n;m_1,...,m_n) \otimes \Delta(l_i,m_1,...,m_n)
$$
So we conclude that (\ref{a0}) is equal to 
$$
- \sum_{i,j =0}^n (-1)^{i+j}  
(l_i\vert l_0,...,\hat 
l_i,...,l_n;m_0,...,\hat m_j,...,m_n) \otimes \Delta(l_i,m_0,...,\hat
m_j,...,m_n)
$$
The considerations for $\nu_{1,n-1}$ are similar. The proposition is proved.

\section {Main construction}

{\bf 1. The weight two case}. Let $B_2(F)$ be the quotient of the free abelian group ${\Z}[F^* \backslash \{1\}]$ generated by the symbols $\{x\}$ where $x \in F^* \backslash \{1\}$ modulo the subgroup $R_2(F)$ generated by the "five term relations", i.e. by the elements
$$
\sum_{i=1}^5 (-1)^i \{r(x_1,...,\hat x_i,...,x_5)\}, \qquad x_i \in {\Bbb P}^1(F), \quad x_i \not = x_j
$$
Denote by $\{x\}_2$ the image of the generator $\{x\}$ in $B_2(F)$. One can prove (see for example [G1]) that there is a well defined homomorphism
$$
\delta_2: B_2(F)  \longrightarrow    \Lambda^2 F^{\ast}, \quad \{x\}_2 \longmapsto (1-x) \wedge x
$$
The complex we get is called the Bloch complex.

In [BGSV1-2] there was defined a  homomorphism of complexes
$$
\begin{array}{ccc} \label {a2}
A_2&\stackrel{\nu}{\longrightarrow} &A_1 \otimes A_1 \\
\downarrow a_2 &&\downarrow a_1 \wedge a_1\\
B_2 &\stackrel{\delta_2}{\longrightarrow} & \Lambda^2 F^{\ast}\\
\end{array}
$$
Namely  
\begin {equation} \label{a51}
a_2(l_0,l_1,l_2;m_0,m_1,m_2) := <l_0,l_1,l_2;m_0,m_1,m_2>_2 :=
\end {equation}
$$
\sum_{i,j=0}^{2}(-1)^{i+j} 
\{r(l_i|l_0,..., \hat l_i,...,l_2;m_0,... ,\hat m_j,...,m_2)\}_2 \in B_2(F)  $$

{\bf 2. The weight three motivic complex related to the classical trilogarithm}. Let $V_3$ be a three dimensional vector space over a field $F$. Choose a volume form $\omega_3 \in det V_3^*$. Recall that for any three vectors $l_1,l_2,l_3$ in $V_3$  we have defined the "determinant" 
$$
\Delta(l_1,l_2,l_3):= <l_1\wedge l_2\wedge l_3, \omega_3 >
$$
 Let us define the generalized cross-ratio of six generic points $x_1,...,x_6$ on the plane $ {\Bbb P}^2$ by setting
\begin{equation} \label{dddd}
r_3(x_1,...,x_6):= \frac{1}{15}{\rm Alt}_6\{\frac{\Delta(\tilde x_1,\tilde x_2,\tilde x_4)
\Delta(\tilde x_2,\tilde x_3,\tilde x_5)\Delta(\tilde x_3,\tilde x_1,\tilde x_6)}{\Delta(\tilde x_1,\tilde x_2,\tilde x_5)\Delta(\tilde x_2,\tilde x_3,\tilde x_6)\Delta(\tilde x_3,\tilde x_1,\tilde x_4)}\} \in {\Z}[F^* \backslash \{1\}]
\end{equation}
Here $\tilde x$ is a vector projecting to the point $x$. The  ratio does not depend on the choice of these vectors. 

{\it The group  $R_3(F)$ of functional equations for the trilogarithm}. I will use the following definition, which is a bit ad hoc. $R_3(F)$ is the group generated by $\{x\}_3 -\{x^{-1}\}_3$, the "seven term" relations (containing actually $7!$ terms)
\begin{equation} \label{mystery}
\sum_{i=1}^7 (-1)^i  r_3(x_1,...,\hat x_i,...,x_7), \qquad x_i \in {\Bbb P}^2(F)
\end{equation}
where the points $x_1,...,x_7$ are in {\it generic} position in the plane, and Kummer's functional equation for the trilogarithm:
$$
K(x,y):=  -\{\frac{x(1-y)^2}{y(1-x)^2}\} - 
 \{xy\}- 
 \{\frac{x}{y}\}- 
2\{1\} +
$$
$$ 2\Bigl(\{\frac{-x(1-y)}{(1-x)}\} +
\{\frac{x(1-y)}{y(1-x)}\} + 
\{\frac{-y(1-x)}{1-y}\} + 
\{\frac{1-x}{1-y}\}+ \{y\} +\{x\}\Bigl)
$$
It might  be true that Kummer's relation  follows from  the generic seven term 
relations  and the one $\{x\}_3 -\{x^{-1}\}_3$. But I just add them 
to the list.

{\bf Remark}. A more natural way to define the group  $R_3(F)$ is this. 
We extend   the generalized cross-ratio to arbitrary configurations of $6$ points on the plane. Then
$R_3(F)$ is given by the seven term relations for arbitrary configurations
of seven points in the plane. Using the main results of [G0-G1] one can show 
that both definitions lead to the same group. In particular we get 
Kummer's relations for a certain degenerate configuration of seven points.
 To extend the definition of $r_3$ we  take the ``limit value'' of the definition (\ref{dddd}) using $\{0\}_3 = \{\infty\}_3=0$. If two of the points $x_i$ coincide or four of them are  on a line then $r_3(x_1,...,x_6)=0$. For the remaining
cases  see lemma \ref{3.7'} below. To make the exposition shorter 
 I will not use this approach. 

Set
$$
B_3(F) := \frac{{\Z}[F^* \backslash \{1\}]}{R_3(F)}
$$

Denote by $\{x\}_3$ the image  in $B_3(F)$ of the generator $\{x\}$. One can prove ([G1-3]) that there is a well defined homomorphism
$$
\delta_3: B_3(F)  \longrightarrow    B_2(F)\otimes  F^{\ast}, \quad \{x\}_3 \longmapsto \{x\}_2 \otimes x
$$
We get the weight three motivic complex related to the classical trilogarithm
$$
B_3(F)  \longrightarrow    B_2(F)\otimes  F^{\ast}   \longrightarrow    \Lambda^3 F^{\ast}
$$
where $\{x\}_2 \otimes y\longmapsto (1-x) \wedge x \wedge y$. 
 
\begin{lemma} \label{3.333} Using the notations introduced in (6)
\begin{equation} \label{33.33}
\delta_3\circ r_3(l_1,l_2,l_3,l_4,l_5,l_6)  =    -\frac{1}{18} {\rm Alt}_6 \Bigl(<l_1,l_2,l_3;l_4,l_5,l_6>_2 \otimes 
\Delta(l_1,l_2,l_3)\Bigr)  =
\end{equation}
$$
-2 <l_1,l_2,l_3;l_4,l_5,l_6>_2\otimes\Delta(l_1,l_2,l_3) + \quad \mbox {$19$ other terms}
$$
\end{lemma}

{\bf Proof}. Let $C_n(V_3)$ be the free abelian group generated by the configurations of $n$ generic vectors  (i.e. $n$-tuples of  vectors in generic position modulo the action of $GL(V_3)$) in $V_3$. Let 
$$
d: C_n(V_3)\to C_{n-1}(V_3), \quad d(l_1,...,l_n) := \sum_{i=1}^n (-1)^{i-1}(l_1,...,\hat l_i,...,l_n)
$$ 
Define a homomorphism $C_5(V_3) \to B_2(F) \otimes F^*$ by setting
$$
f_5(3)(l_1,...,l_5) := \frac{1}{2} {\rm Alt}_5\Bigl( \{r(l_1|l_2,...,l_5)\}_2 \otimes \Delta(l_1,l_2,l_3)\Bigr)
$$
According to theorem 2.3 in the Appendix to [G3]  the following diagram is commutative:
$$
\begin{array}{ccc}
C_6(V_3)&\stackrel{d}{\longrightarrow}&C_5(V_3)\\
&&\\
r_3 \downarrow& &\downarrow f_5(3)\\
&&\\
B_3(F)&\stackrel{\delta_3}{\longrightarrow}&B_2(F) \otimes F^*
\end{array}
$$
Therefore 
$$
\delta_3\circ r_3(l_1,l_2,l_3,l_4,l_5,l_6)  = -\frac{1}{2} {\rm Alt}_6 \Bigl(\{r(l_1|l_2,l_3,l_4,l_5)\}_2 \otimes 
\Delta(l_1,l_2,l_3)\Bigr)
$$
which is equal  to (\ref{33.33}).

{\bf 3.  A homomorphism  between the weight three motivic complexes}. We have constructed in chapter 2 a complex
$$
A^0_3 \stackrel{\nu_{2,1} \oplus  \nu_{1,2}}{\longrightarrow} A_2 \otimes A_1 \oplus A_1
\otimes A_2 \stackrel{ \nu_{1,1} \otimes Id - Id \otimes
\nu_{1,1}}{\longrightarrow} A_1 \otimes A_1 \otimes A_1
$$
Let us  define a homomorphism  of complexes
$$
\begin{array}{ccccc} \label {a3}
A^0_3&\stackrel{\nu}{\longrightarrow} &A_2 \otimes A_1 \oplus A_1
\otimes A_2 & {\longrightarrow}& A_1 \otimes A_1 \otimes A_1\\
&&&&\\
\downarrow a_3 &&\downarrow a_2 \wedge a_1 &&\downarrow \wedge^3 a_1\\
&&&&\\
B_3 &\stackrel{\delta_3}{\longrightarrow} & B_2(F) \otimes F^{\ast}&
\stackrel{ }{\longrightarrow} & \Lambda^3 F^{\ast}\\
\end{array}
$$
where $a_2 \wedge a_1(x_2 \otimes x_1 + y_1 \otimes y_2):= a_2(x_2) \otimes
a_1(x_1) - a_2(y_2) \otimes  a_1(y_1)$ and $\wedge^3a_1(x_1 \otimes x_2 \otimes x_3) = a_1(x_1) \wedge a_1(x_2) \wedge a_1(x_3)$.  

\begin{center}
\hspace{4.0cm}
\epsffile{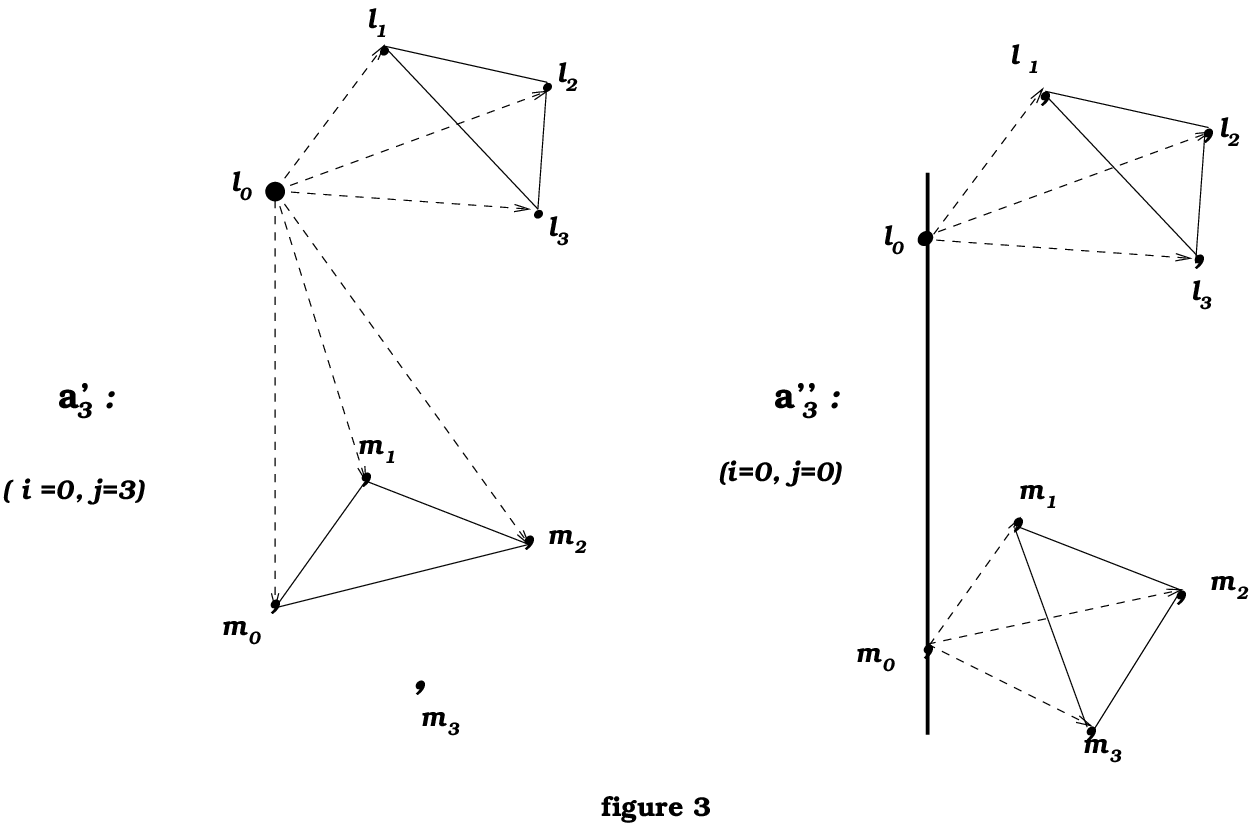}
\end{center}

Let 
$x_1,x_2,x_3;y_1,y_2,y_3$ be six different points on a line. Set

$$
 \mu_3 (x_1,x_2,x_3;y_1,y_2,y_3) = 
$$
$$
\frac{1}{4}\cdot  
{\rm Alt}_{(x_1,x_2,x_3)(y_1,y_2,y_3)}\Bigl(
\{r(x_1,y_2,x_2,y_1)\}_3 - \{r(x_1,y_1,x_2,y_2)\}_3\Bigr)
$$

This formula, which contains 18 terms,  means  simply that $ \mu_3
(x_1,x_2,x_3;y_1,y_2,y_3)$ 
is skewsymmetric with respect to $x_1,x_2,x_3$ and $y_1,y_2,y_3$ and a typical
term is $\{r(x_1,y_2,x_2,y_1)\}_3$.
 
We set $$
a_3 := \frac{1}{6} a_3' -\frac{1}{3}a_3''
$$
 where $a_3'$ and $a_3''$ are
defined on generators by the following formula (see fig. 3):
$$
a_3' (l_0,l_1,l_2,l_3;m_0,m_1,m_2,m_3) : =  \sum_{i,j=0}^{3}(-1)^{i+j} 
r_3(l_i|l_0,..., \hat l_i,...,l_3;m_0,..., \hat m_j,...,m_3)   
$$
$$
a_3'' (l_0,...,l_3;m_0,...,m_3) : =  \sum_{i,j=0}^{3}(-1)^{i+j}  \cdot 
   \mu_3(l_i,m_j|l_0,..., \hat l_i,...,l_3;m_0,..., \hat m_j,...,m_3) 
$$

Denote  by $\tilde A^0_3$ the free abelian group generated by the generators of the group
$ A^0_3$, i.e. by pairs of tetrahedra in generic position.

\begin {theorem} \label{Di}
The diagram
$$
\begin{array}{ccc} \label {a31}
\tilde A^0_3&\stackrel{\nu}{\longrightarrow} &A^0_2 \otimes A^0_1 \oplus A^0_1
\otimes A^0_2 \\
&&\\
\downarrow a_3 &&\downarrow a_2 \wedge a_1 \\
&&\\
B_3 &\stackrel{\delta_3}{\longrightarrow} & B_2(F) \otimes F^{\ast}
\end{array}
$$
is commutative.
\end {theorem}

{\bf Remark}. The proof of this 
  theorem is the crucial point of the paper. The commutativity  of the diagram is what we really wanted  for the homomorphism $a_3$. Surprisingly neither  homomorphism $a_3'$ nor $a_3''$   make the diagram commutative, even up to a scalar. Only their sum does the job. Moreover,  theorem (\ref{Di})
"morally" implies that $a_3$ should be a homomorphism of groups. 
We will prove this later on.

{\bf Proof}. Because of the skew-symmetry of the
formulas   for $\nu$ (see proposition  \ref{cop} ) it is sufficient to compute for 
$$
(a_2 \wedge a_1)\circ \nu(l_0, ..., l_3;m_0,...,m_3) \quad \mbox{ and } \quad  \delta_3 \circ
a_3(l_0, ..., l_3;m_0,...,m_3)
$$ 
the $B_2(F)$-factors  of the following elements of $F^*$:
$$
\Delta(m_0,m_1,m_2,m_{3}), \quad \Delta(l_0,m_0,m_1,m_2), \quad \Delta(l_0, l_1,m_0, m_{1 }), 
$$
$$
\quad \Delta(l_0,l_1,l_2,m_{0}),\quad \Delta(l_0,l_1,l_2,l_{3})
$$  
 Since $(a_2\wedge a_1) \circ \nu(L,M) = (a_2\wedge a_1) \circ  
\nu (M,L)$ by the definition of $a_2$ and $a_1$ and 
by proposition 2.3,     
and $a_3(L;M) = a_3(M;L)$ by lemma (\ref{A}),  we see that it is
sufficient to consider only the  first three of them. 

 It follows from proposition (\ref{cop}) that in $(a_2 \wedge a_1)\circ
\nu(l_0, ..., l_3;m_0,...,m_3)$ 
appears only  
\begin{equation} \label{**<>}
<l_0|l_1,l_2,l_3;m_1,m_2,m_3>_2 \otimes \Delta(l_0,m_1,m_2,m_3)
\end{equation}

{\it Step} 1. Let us do the computations for $\delta_3 \circ a_3(l_0, ...,
l_3;m_0,...,m_3)$. The crucial and most nontrivial case is the term  with
$\Delta(l_0,l_1,m_0,m_1)$. If the diagram is commutative  it must be zero because of the observation we just made. 
The only summands in 
 $a_3'(l_0, ..., l_3;m_0,...,m_3)$ that give a contribution to this term are:
$$
- r_3(l_0|l_1,l_2,l_3,m_0,m_1,m_2) + r_3(l_0|l_1,l_2,l_3,m_0,m_1,m_3) +
$$
$$
r_3(l_1|l_0,l_2,l_3,m_0,m_1,m_2) - r_3(l_1|l_0,l_2,l_3,m_0,m_1,m_3) 
$$

  Lemma (\ref{3.333}) shows using the symmetry $1<->4, 2<->5$ that the term  with
$\Delta(l_0,l_1,m_0,m_1)$  in $\delta_3 \circ a_3(l_0, ...,
l_3;m_0,...,m_3)$  is
$$
-2 \Bigl(-<l_0|l_1,m_0,m_1;l_2,l_3,m_2>_2  + <l_0|l_1,m_0,m_1;l_2,l_3,m_3>_2+
$$
$$
<l_1|l_0,m_0,m_1;l_2,l_3,m_2>_2 - <l_1|l_0,m_0,m_1;l_2,l_3,m_3>_2\Bigr)
\otimes \Delta(l_0,l_1,m_0,m_1) = 
$$
\begin {equation} \label{c41}
2 \cdot {\rm Alt}_{l_0l_1}{\rm Alt}_{m_2m_3}\Bigl( <l_0|l_1,m_0,m_1;l_2,l_3,m_2>_2\Bigr)\otimes \Delta(l_0,l_1,m_0,m_1)
\end  {equation}

{\it Step} 2. The contribution for $\Delta(l_0,l_1,m_0,m_1)$  in $\delta_3\circ a_3''(l_0, ..., l_3;m_0,...,m_3)$ is coming from 
$$
\delta_3\circ  \mu_3\Bigl( (l_0m_0|l_1,l_2,l_3;m_1,m_2,m_3) - (l_0m_1|l_1,l_2,l_3;m_0,m_2,m_3)
$$
$$
- (l_1m_0|l_0,l_2,l_3;m_1,m_2,m_3)  + (l_1m_1|l_0,l_2,l_3;m_0,m_2,m_3)\Bigr) 
$$
And the contribution is equal to 
$$
{\rm Alt}_{(l_0,l_1);(m_2,m_3);(m_0,m_1);(l_2,l_3)}\Bigl( (l_0m_0|l_1,m_1,l_2,m_2)_2  
\Bigr)\otimes \Delta(l_0,l_1,m_0,m_1) = 
$$
$$
{\rm Alt}_{(l_0,l_1);(m_2,m_3) }\Bigl( <l_0|m_0,l_1,m_1;l_2,m_2,l_3>_2  
\Bigr)\otimes \Delta(l_0,l_1,m_0,m_1) =
$$
$$
{\rm Alt}_{(l_0,l_1);(m_2,m_3) }\Bigl( <l_0|l_1,m_0, m_1;l_2, l_3,m_2>_2  
\Bigr)\otimes \Delta(l_0,l_1,m_0,m_1)
$$
 (We use a shorthand $(l_0m_0|l_1,m_1,l_2,m_2)_2$ for $\{r(l_0m_0|l_1,m_1,l_2,m_2)\}_2$). Comparing the last formula with (\ref{c41}) we see that 
$$
\mbox{{\it    $\frac{1}{2}a_3' - a_3''$ has contribution $0\otimes \Delta(l_0,l_1,m_0,m_1)$}}
$$

{\it Step} 3. Now look at the terms of $\delta_3
 \circ
a_3(l_0, ..., l_3;m_0,...,m_3)$ with $\Delta( m_0, m_1, m_2,  m_{3 })$ and $\Delta( l_0, m_1, m_2, m_3)$. 

1) $\Delta(m_0,m_1,m_2,m_3)$ does not appear in $a_3(l_0, ...,
l_3;m_0,...,m_3)$.

2i)Using proposition (\ref{cop})    we get that 
$\Delta(l_0,m_1,m_2,m_3)$ appears in $\delta_3 \circ a_3'(l_0, ..., l_3;m_0,...,m_3)$ as 
$$
2<l_0|l_1,l_2,l_3;m_1,m_2,m_3>_2 \otimes \Delta(l_0,m_1,m_2,m_3)
$$

2ii) Let  us  compute $\delta_3 \circ a_3''(l_0, ..., l_3;m_0,...,m_3)$. 
There are only $3$ terms in $a_3''$ where the term   $\Delta(l_0,m_1,m_2,m_3)$ appears; we can write them as
$$
- 1/2{\rm Alt}_{(m_1,m_2,m_3)}\mu_3(l_0m_1|l_1,l_2,l_3;m_0,m_2,m_3)
$$
It is equal to 
$$
  1/2 {\rm Alt}_{(l_1,l_2,l_3)(m_1,m_2,m_3)} (l_0m_1|l_1,m_2,l_2,m_3)_2\otimes
\Delta(l_0,m_1,m_2,m_3) = 
$$
$$
-2 <(l_0 |m_1,m_2,m_3;l_1,l_2,l_3)>_2\otimes\Delta(l_0,m_1,m_2,m_3)
$$
For the last step use $(l_0m_1|l_1,m_2,l_2,m_3)_2 = -(l_0m_1|m_2,m_3,l_1,l_2)_2$. 
So  the contribution  of $\frac{1}{2}a_3'-a_3''$ is
$$
3\cdot  <l_0\vert l_1,l_2,l_3;m_1,m_2,m_3>_2 \otimes \Delta(l_0,m_1,m_2,m_3)
$$
It remains to compare this answer with  (\ref{**<>}).
Theorem (\ref{Di}) is proved.

\begin {theorem} \label{7.1}
$a_3$ is a homomorphism of groups.
\end {theorem}

We start the proof with 

\begin {proposition} 
Both $a_3'$ and $a_3''$ send relations 2), 3), 4) to zero.
\end {proposition}

{\bf Proof}. This is clear for relations 2) and 4). It is also clear
that $a_3''$ is additive in $L$ and in $M$. To check that $a_3''$
sends  
the dual additivity in $L$ to zero notice that the typical term for  
$a_3''(\sum_{i=0}^{3}(-1)^{i}  \cdot(l_i|l_0,...,\hat
l_i,...,l_4;m_0,...,m_3))$ is  
$$
\pm \mu_3(l_il_jm_k|l_0,...,\hat l_i,...,\hat
l_j,...,l_4;m_0,...,\hat m_k,...,m_3)
$$
It is symmetric in $l_i,l_j$. So we get 0 after alternation in
$(l_0,...,l_4)$.

Similarly using the skewsymmetry of $r_3 $ we immediately see that the map   $a_3'$
sends  dual additivity in $L$ to zero. The additivity in
$M$ of $a_3'$ 
is obvious.  
Therefore the proposition  follows immediately from the following lemma.

\begin {lemma} \label{A}
 $a_3'(L;M) = a_3'(M;L); \qquad a_3''(L;M) = a_3''(M;L)$ 
\end {lemma}

{\bf Proof}. The statement about $a_3''$ is clear from the definition. 
Let us prove it for $a_3'$. Applying the 7-term relation for the
configuration of 7 points in ${\Bbb P}^2$
$(l_i|l_0,...,\hat l_i,...,l_3,m_0,...,m_3)$  we get 
\begin  {equation} \label{fedya}
a_3'(l_0,,...,l_3;m_0,...,m_3) = - \sum_{0 \leq i \not = j \leq 3 } \gamma(i,j) 
r_3(l_i|l_0,...,\hat l_i,..., \hat l_j,...,l_3,m_0,...,m_3)
\end {equation}
where $\gamma(i,j) = (-1)^{i+j}$ if $i<j$    and $-(-1)^{i+j}$ otherwise.
Applying the dual 7-term relation to the configuration 
$(l_0,...,\hat l_i,...,l_3,m_0,...,m_3)$  of 7 points in ${\Bbb P}^3$ we can
rewrite (\ref{fedya}) as 
$$
- \sum_{i,j=0}^3 (-1)^{i+j}r_3
(m_i|l_0,...,\hat l_j,...,l_3,m_0, ...,\hat m_i,...,m_3)
$$
It remains to interchange the $l$- and $m$- tetrahedra  using  the
skewsymmetry. 
The lemma is proved.

Neither $a_3'$ nor $a_3''$ send the relation  1) to zero.  Only their weighted sum $a_3$ does the job. To show this we have to prove the proposition 
\ref{3.7}  below.

Let $L_0,...,L_3$ be  4 lines and $m_0,...,m_3$ 4 points in ${\Bbb P}^2$.
Let $$
(L_0,...,\hat L_i,...,L_3;m_0,...,\hat m_j,...,m_3)
$$ be a pair of
triangles where the first triangle is given by its sides and the second
one by its vertices (see fig. 4). For example  
$$
(L_1,L_2,L_3;m_1,m_2,m_3) = (l_{23},l_{13},l_{12},m_1,m_2,m_3)
$$
where $l_{ij}:= L_i \cap L_j$, so 
the first three points are the vertices of the triangle $(L_1,L_2,L_3)$.


\begin {proposition} \label{3.7}
 Let $L_0,...,L_3$ be  4 lines and $m_0,...,m_3$ 4 points in ${\Bbb P}^2$. 
Then 
\begin {equation} \label{aa4}
 \sum_{i,j =0}^3 (-1)^{i+j} 
r_3(L_0,...,\hat L_i,...,L_3;m_0,...,\hat m_j,...,m_3) =
\end  {equation}
\begin {equation} \label{aa54}  
2 \sum_{i,j =0}^3 (-1)^{i+j} 
    \mu_3(m_i|l_{j0},...,\hat l_{jj},...,l_{j3};m_0,...,\hat m_i,...,m_3)
\end  {equation}
\end {proposition}

\begin{center}
\hspace{4.0cm}
\epsffile{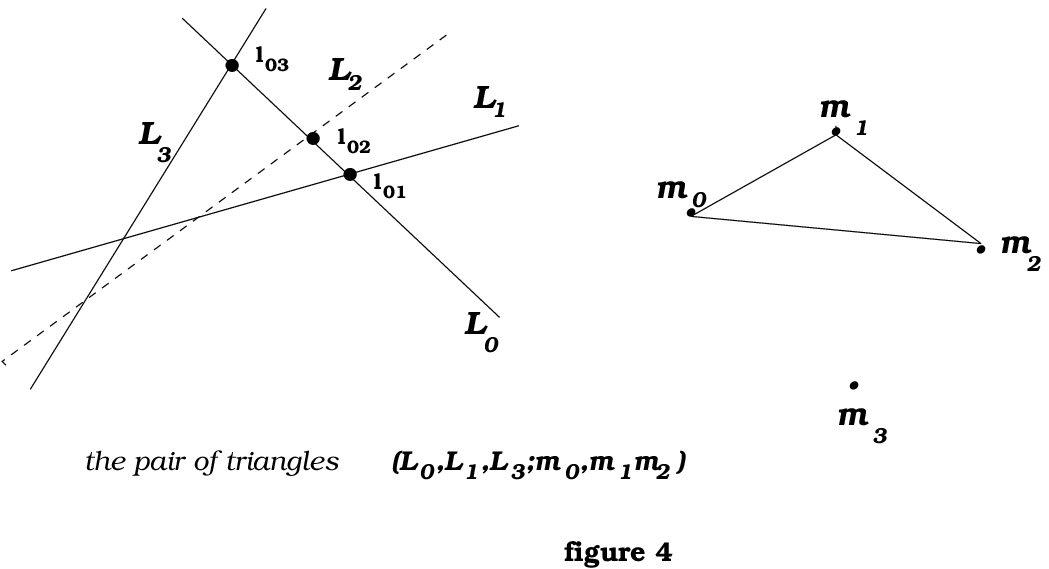}
\end{center}

{\bf Proof}.
   Applying the 7-term relation to the configuration of 7
points formed by the
 vertices of the triangle $L_0,...,\hat L_j,...,L_3$ and $m_0,...  ,m_3$ we  rewrite (\ref{aa4}) as  
$$
-\frac{1}{2}{\rm Alt}_{(L_0,...,L_3)}r_3(l_{13},l_{23},m_0,m_1,m_2,m_3)   
$$
(We alternate the indices $0,1,2,3$ in the $l$-variables only). It is equal to
$$
- r_3(l_{13},l_{23}, m_0,m_1,m_2,m_3) +  
   r_3( l_{03},l_{23} ,m_0,m_1,m_2,m_3) 
-  r_3(l_{03},l_{13}, m_0,m_1,m_2,m_3)   
$$
$$
\quad + \mbox{ 9 other terms}
$$
Consider the three intersection points  of the line $L_i$ with  the
other lines and add to them the points $m_0,m_1,m_2,m_3$. For instance $(l_{03},l_{13},l_{23},m_0,m_1,m_2,m_3)$ is the configuration related to the line $L_3$. Applying the 7-term relation to these
  configurations  we rewrite the previous formula as  
\begin {equation} \label{aaa}
 - r_3\Bigl(\sum_{i=0}^3 (-1)^i (l_{03},l_{13},l_{23},m_0,...,\hat m_i,...,m_3) - \sum_{i=0}^3 (-1)^i (l_{10},l_{20},l_{30},m_0,...,\hat m_i,...,m_3)
\end {equation}
 \begin {equation} \label{aab}
 +\sum_{i=0}^3 (-1)^i (l_{21},l_{31},l_{01},m_0,...,\hat m_i,...,m_3)
 -\sum_{i=0}^3 (-1)^i (l_{32},l_{02},l_{12},m_0,...,\hat m_i,...,m_3)\Bigr)
 \end {equation}
  
Now comes a little trick: we will use the fact that one can extend the generalized cross ratio to certain degenerate configurations of six points such that the seven term relation holds. 
\begin {lemma} \label{3.7'}
 Suppose the points $l_1,l_2,l_3$ are on the same line. Then
$$
r_3(l_1,l_2,l_3,m_1,m_2,m_3) = - 2 \mu_3(l_1,l_2,l_3;m_1,m_2,m_3)
$$
\end {lemma}
 
{\bf Remark}. It is easy to check that
$$
\delta_3 \circ (r_3 + 2\mu_3) (l_1,l_2,l_3,m_1,m_2,m_3) =0
$$
This implies that $(r_3 + 2\mu_3) (l_1,l_2,l_3,m_1,m_2,m_3)$ is a functional equation for the trilogarithm. The following proof verifies that it is already in the group $R_3$ defined above  (so we do not have to enlarge this group).
 
{\bf Proof}.  The following proof based on  direct computation using  $\{0\}_3 = \{\infty\}_3 =0$ was provided to me by Herbert Gangl.
Consider the following special configuration of 6 points on the plane (fig. 5)
\begin{center}
\hspace{4.0cm}
\epsffile{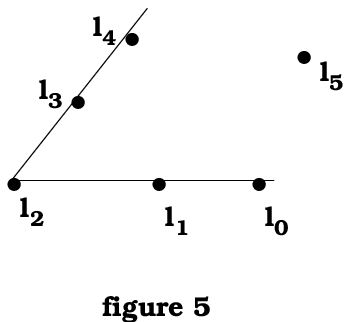}
\end{center}
given by the columns of the matrix
$$ {\cal C}(b,c):=
\left (\matrix{0&0&0&1&1&1 \cr 
1& 1&0&1&1&0 \cr
 0&c&
    1&b&1&0\cr}\right ) 
$$

It is sufficient to prove the lemma for this configuration, since the general case will follow by the seven term relation.

Set 
$$
\tilde l_3\{x\}:=  \{1-x^{-1}\}_3 - \{1-x\}_3
$$
$$
M_3({\cal C}(b,c)):= \tilde l_3\Bigl(\{\frac{b-c}{1-c}\}_3 - \{b\}_3+ \{\frac{b(1-c)}{b-c}\}_3
- \{1-c\}_3 + \{\frac{b-c}{b}\}_3\Bigr)
$$

A configuration of six points  $x_1,x_2,x_3,l_1,l_2,l_3$ on the plane, three 
of which, $x_1,x_2,x_3$,  are on a line $L$, is determined completely by the 
configuration \linebreak $(x_1,x_2,x_3,n_1,n_2,n_3)$ of $3+3$ points on the line $L$,  where 
the point $n_1:= l_2l_3 \cap L$ and so on, 
see fig. 6.

\begin{center}
\hspace{4.0cm}
\epsffile{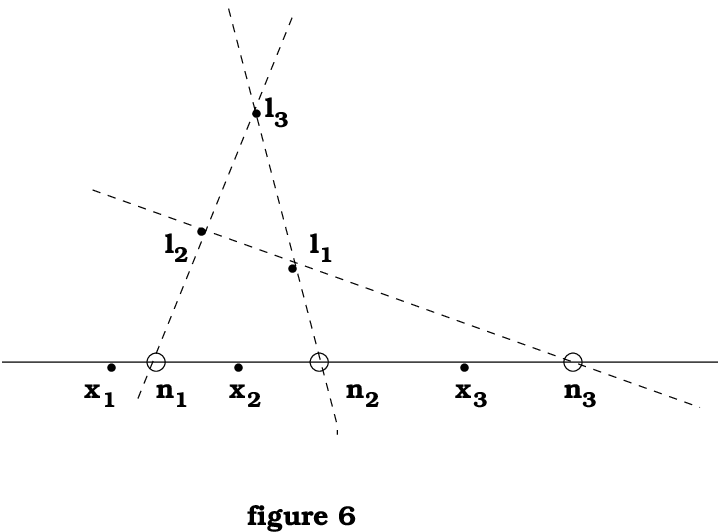}
\end{center}
(Recall that by configuration we always mean its projective equivalence class.) 

Let $\tilde {\cal C}(b,c)$ be the configuration of six points on the line 
corresponding to the configuration ${\cal C}(b,c)$ by this rule. Then 
$$
M_3({\cal C}(b,c)) = -\mu_3(\tilde {\cal C}(b,c))
$$

One has
$$
15 \cdot r_3({\cal C}(b,c)) = -18 \{1-b\}_3 -12\{\frac{b}{1-c}\}_3 -18 
\{\frac{b-c}{b-1}\}_3 - 18 \{c\}_3 
$$
$$
- 12\{\frac{1-b}{c}\}_3 - 24\{\frac{b(c-1)}{c(b-1)}\}_3 - 6 \{\frac{c}{(c-b)(1-b)}\}_3 - 6 \{\frac{c(b-c)}{ b-1 }\}_3 
- 18 \{\frac{c-b}{c}\}_3
$$
$$
+18 \{1-b^{-1}\}_3 + 12\{b-c\}_3 +18\{\frac{b-1}{c-1}\}_3  +12\{\frac{b-1}{c}\}_3 
+18\{\frac{b}{c}\}_3 
$$
$$
+24\{\frac{b-c}{c(b-1)}\}_3 +6\{\frac{b(b-1)}{c(c-1)}\}_3
 +18\{1-c^{-1}\}_3 +6\{\frac{(1-b)(1-c)}{b\cdot c}\}_3
$$
Recall two Kummer  relations:
$$
6 K(\frac{b}{1-c},\frac{c}{1-b}) = 12\{\frac{b}{1-c}\}_3 + 12\{\frac{c}{1-b}\}_3 + 
  12\{\frac{1-b}{1-c}\}_3 + 12\{\frac{c}{c-1}\}_3 + 12\{\frac{b}{b-1}\}_3 
$$
$$
+ 12\{\frac{c}{b}\}_3 - 6\{\frac{bc}{(1-b)(1-c)}\}_3 -6\{\frac{b(b-1)}{c(c-1)}\}_3 -
6\{\frac{c(1-b)}{b(1-c)}\}_3 -12\{1\}_3
$$
and
$$
-6 K(b-c,\frac{b-1}{c}) = -12\{b-c\}_3 - 12\{\frac{b-1}{c}\}_3
  -12\{c\}_3 - 12\{1-b\}_3 - 12\{\frac{c}{c-b}\}_3 
$$
$$
- 12\{\frac{b-1}{b-c}\}_3 + 6\{\frac{(c-b)(1-b)}{c}\}_3 +6\{\frac{c(b-c)}{b-1}\}_3 +
6\{\frac{c(b-1)}{b-c}\}_3 +12\{1\}_3
$$
Then one computes, adding the three expressions above,
$$
15r_3({\cal C}(b,c)) + 6K(\frac{b}{1-c},\frac{c}{1-b})-6 K(b-c,\frac{b-1}{c}) = 
$$
$$
30\Bigl( -\{1-b\}_3 + \{\frac{b}{b-1}\}_3 - \{\frac{b-c}{b-1}\}_3 +
 \{\frac{1-b}{1-c}\}_3 
 - \{\frac{c(1-b)}{b(1-c)}\}_3  + \{\frac{c-b}{c(1-b)}\}_3 
$$
$$
-\{c\}_3 +\{1-c^{-1}\}_3 -\{\frac{c-b}{c}\}_3 +\{\frac{c}{b}\}_3\Bigr) = -30 \mu_3
({\cal C}(b,c))
$$
which is exactly what we wanted. The lemma is proved. Comparing this lemma with (\ref{aaa}), (\ref{aab}) we get proposition 3.4 and hence  
 theorem (\ref{7.1}).
  
$$
$$
Another way to proceed. The right hand side of (\ref{aa4}) considered as a
function of 4 lines and 4 points on the plane satisfies the 5 term
relation with respect to lines as well as points. Using this 
observation it is sufficient to check the proposition for  the
degenerate configuration of 
lines and points when  3 of the lines $L_i$ pass through a point and 3
points among the $m_i$ are on 
a line (see fig. 7).  

\begin{center}
\hspace{4.0cm}
\epsffile{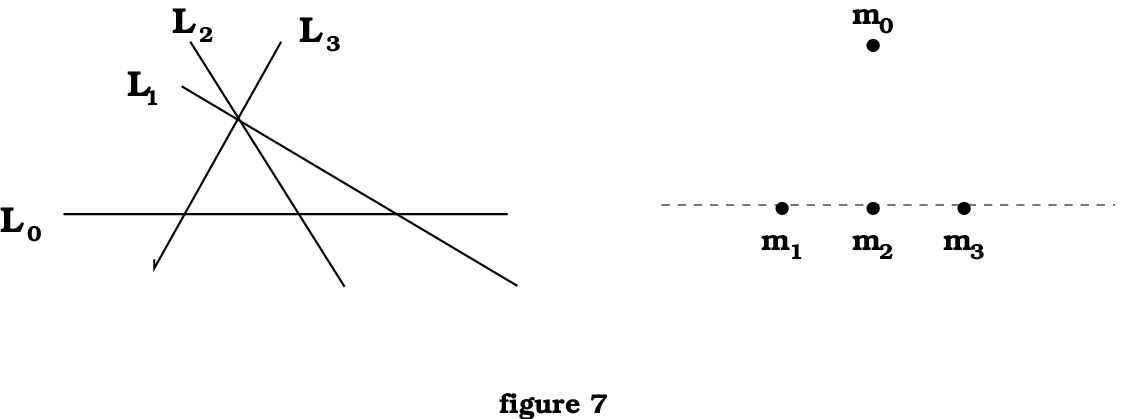}
\end{center}

{\bf 4. Another way to define the homomorphism $a_3$}. Let us associate to $8$ points $(l_0,... ,l_3;m_0,... ,m_3)$ in ${\Bbb P}^3$   a degenerate configuration of $6$ points on the plane, denoted
$$
d(l_i,m_j||l_0,..., \hat l_i,...,l_3;m_0,...,\hat m_j,...,m_3)
$$
as follows. Consider the configuration
$$
(l_i|l_0,..., \hat l_i,...,l_3;m_0,...,m_3)
$$
of $7$ points in ${\Bbb P}^2$. We construct out of them a degenerate 
configuration of $6$ lines in the plane. Let $L_{ij}$ (resp ($M_{ij}$) be the line through the points $l_i$ and
$l_j$ (resp. $m_i$ and
$m_j$). Take the three lines 
formed by the sides of the triangle $(l_i|l_0,..., \hat l_i,...,l_3)$ and 
add to them the three lines $M_{j0},...\hat M_{jj},...,M_{j3}$ connecting $m_j$ with 
the other three $m$-points.  For example
$$
d(l_0,m_0||l_1,l_2,l_3;m_1,m_2,m_3) =
 (L_{23},L_{13},L_{12},M_{01},M_{02},M_{03})
$$
Let $\bar C_{2n}(F)$ be the free abelian group generated by the configurations of 
{\it arbitrary} $2n$ $F$-points in ${\Bbb P}^{n-1}$.
One can define $a_3$ as a composition
$$
A^0_3(F) \stackrel{p_3}{\longrightarrow} \bar C_6(F)/\{\mbox{7 term relations}\}
\stackrel{r_3}{\longrightarrow} B_3(F) 
$$
where (see fig. 8)
$$
p_3(l_0,l_1,l_2,l_3;m_0,m_1,m_2,m_3):= \quad
 \sum_{i,j=0}^{3}\Bigl((-1)^{i+j} 
(l_i\vert l_0,..., \hat l_i,...,l_3;m_0,..., \hat m_j,...,m_3) 
$$
\begin{equation} \label{12}
+ 2 \cdot 
d(l_i,m_j||l_0,..., \hat l_i,...,l_3;m_0,...,\hat m_j,...,m_3)\Bigr)
\end{equation}

\begin{center}
\hspace{4.0cm}
\epsffile{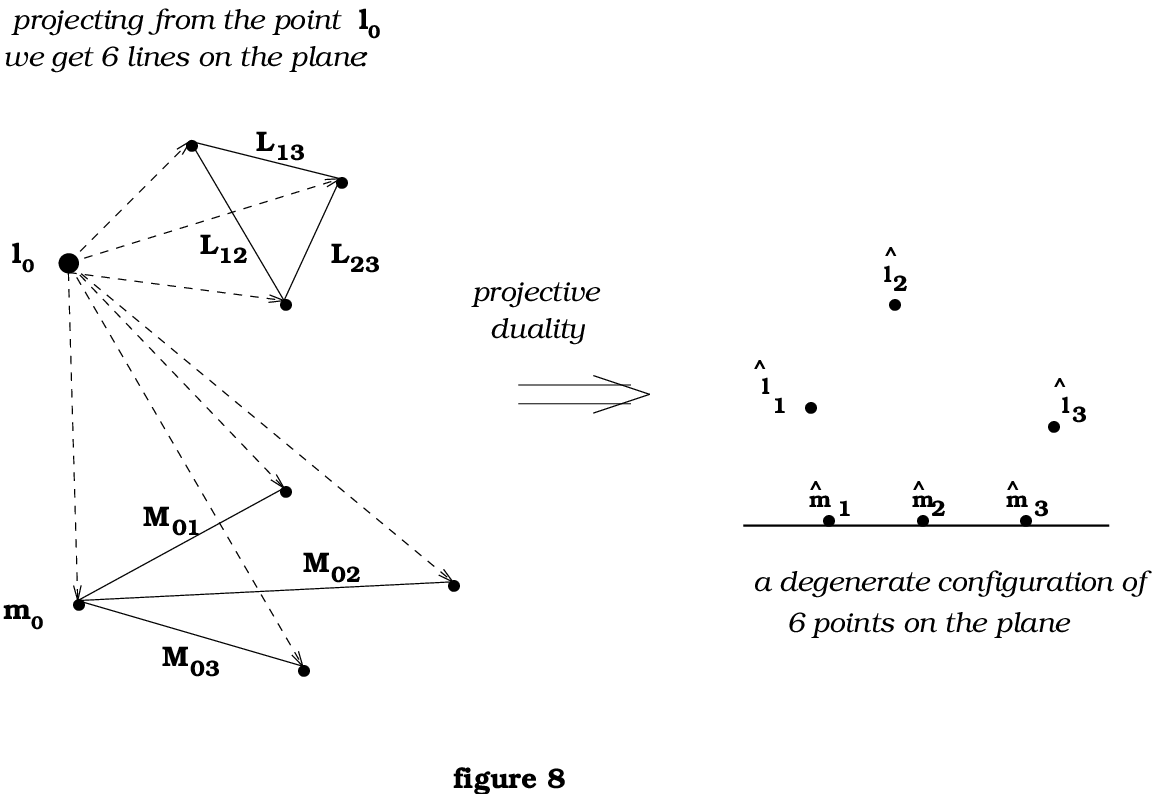}
\end{center}

There is a different candidate for the term (\ref{12}). We define a configuration
\begin{equation} \label{13}
\hat d(l_i,m_j||l_0,..., \hat l_i,...,l_3;m_0,...,\hat m_j,...,m_3)
\end{equation}
of $3+3$ points on a line as follows. Consider the planes $L_i$ and $M_j$ in ${\Bbb P}^3$ (see fig. 9 for $i=0, j=0$.)  Their intersection is a line $L_i \cap M_j$. The sides of the triangles $(l_0,..., \hat l_i,...,l_3)$ and $(m_0,...,\hat m_j,...,m_3)$ cut this line in $3+3$ points. This is the configuration (\ref{13}) we promised to define.

\begin{center}
\hspace{4.0cm}
\epsffile{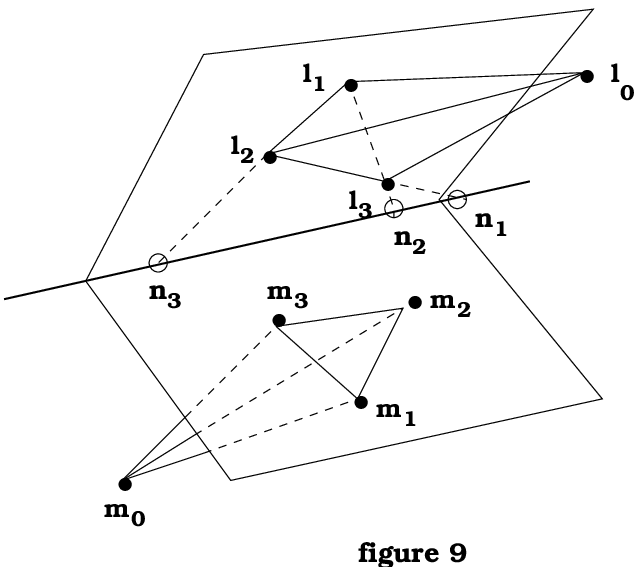}
\end{center}

Recall that we can think of a configuration of $3+3$ points on a line as of a configuration of six points on a plane (see fig. 6). 
So we can describe configuration (\ref{13}) by the configuration of 
points  $(n_0, ... , \hat n_i, ...,n_3, m_0,..., \hat m_j,...,m_3)$ on the plane $M_j$, or by a similar configuration on the plane $L_i$.

The definition of (\ref{13}) given on fig. 9 is 
projectively dual to the one on fig. 8. This means that if we consider the points $l_i, m_j$ as planes in the dual space $\hat {\Bbb P}^3$ then (\ref{12}) corresponds to 
(\ref{13}). 

{\bf Remark}. One has (mainly thanks to lemma (\ref{3.7'}))
\begin{equation} \label{14}
 \sum_{i,j=0}^{3} (-1)^{i+j} \Bigl(
d(l_i,m_j||l_0,..., \hat l_i,...,l_3;m_0,...,\hat m_j,...,m_3) -
\end{equation}
$$
\hat d(l_i,m_j||l_0,..., \hat l_i,...,l_3;m_0,...,\hat m_j,...,m_3)
\Bigr) =0
$$
Let us imagine that we would be able to prove this using only the (possibly degenerate)  seven term relations but not 
 lemma (\ref{3.7'}).  Then it is straightforward to show that $p_3$ sends to zero 
all the defining relations for the group $A_3$ {\it except} 
the degeneracy relation.  So we can define the group $R_3$ by adding to the 
(possibly degenerate) seven term relations the image of the degeneracy relations,
 i.e. the relations  from  proposition (\ref{3.7}).
Moreover we get a nice free gift: now we can  skip lemma (\ref{3.7'}) together 
with its  computational proof. So
 instead of the mysterious relation (\ref{mystery}) for the classical
 trilogarithm which we think of as a relation 
$$
\mbox{generic configuration of $6$ points on ${\Bbb P}^2$ } \quad = \sum \quad \mbox{of degenerate configurations}
$$
 (and which {\it does not} follow from the seven term relations between 
the configurations of six points on the plane), we would have  the 
geometrically natural relations from  proposition (\ref{3.7}). 
However I was not able to prove (\ref{14}) without using lemma (\ref{3.7'}). 
I wish somebody will try. 

{\bf 5. The key ansatz from [G0-3]}. 
Consider the following admissible pair of tetrahedra which represents $0$ on $A_3$ (see fig 10). 

\begin{center}
\hspace{4.0cm}
\epsffile{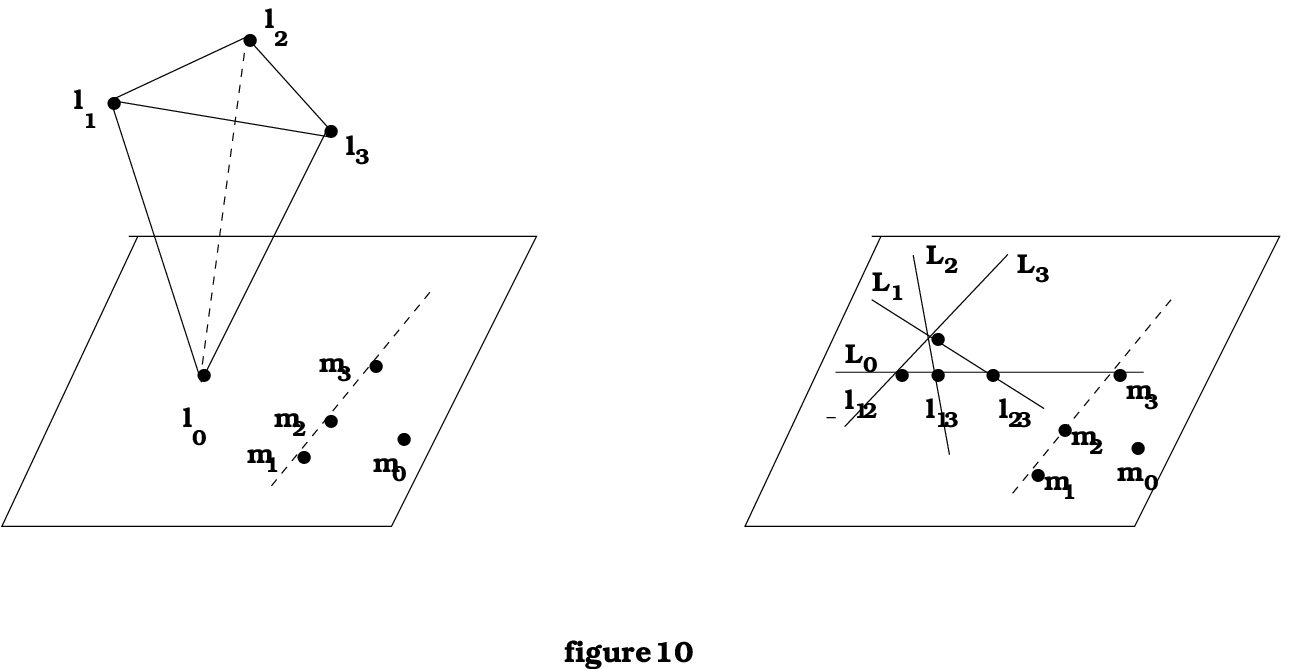}
\end{center}

Here $l_0, m_0, m_1, m_2, m_3$ are on the same plane. The right hand side of the picture illustrates the computation of the homomorphism $a_3$.

Then 
$$
a_3'(L,M) = \{r_3(l_{12}, l_{13}, l_{23}, m_1, m_2, m_3)\}
$$
Indeed, all the other terms in (\ref{aaa})-(\ref{aab}) vanish since by definition 
a degenerate configuration of six points on the plane represents   zero 
in $B_3$ if two of them coincide or four are on the same line. 
The condition that for such $(L,M)$ one has 
$6 a_3(L,M) = (a_3' - 2 a_2'')(L,M) =0$ just means that we get a formula expressing the configuration $l_{12}, l_{13}, l_{23}, m_1, m_2, m_3$ as a sum of configurations like the one on fig. 11 (each of them corresponds to a generator $\{x\}_3$ of the group $B_3$). 

\begin{center}
\hspace{4.0cm}
\epsffile{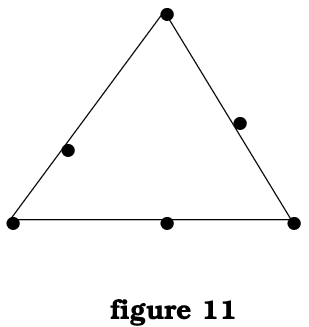}
\end{center}

{\bf 6. A formula for $\Lambda_3(L,M)$}. One can define (see [BGSV]) a commutative product map
$$
\mu: A_k \otimes A_l \to A_{k+l}
$$
Then one should have a structure of a Hopf algebra on $A_{\bullet}$ given by the coproduct $\nu$ and the product $\mu$ (however at the present moment  it is not clear how to define the coproduct for the degenerate admissible pairs of simplices, and thus we can not prove  that $\nu \circ \mu = \mu \otimes \mu (\nu)$).
Let $$
P_n:= \oplus_{k+l=n}\mu(A_k \otimes A_l)\qquad (k>0,l>0)
$$

Consider the homomorphism
$\pi_3: A_3 \to A_3$ given by
$$
\pi_3 (x):= x - \frac{1}{2} \mu \circ \nu (x) + \frac{1}{3}\mu \circ \nu_{1,1,1} (x)
$$
where $\nu_{1,1,1}: A_3 \to A_1 \otimes A_1 \otimes A_1 $ is the composition

$$
A_3 \stackrel{\nu_{2,1}}{\lra} A_2 \otimes A_1 \stackrel{\nu_{1,1}\otimes id}{\lra} A_1 \otimes A_1 \otimes A_1
$$

\begin{lemma}
$\pi_3$ is a projector and $Ker \pi_3 = P_3$. 
\end{lemma}

See also section 4.3 for a generalization.

{\bf Proof}. Let us check that $P_3 \subset Ker \pi_3$. This implies  
$\pi_3^2 = \pi_3$. Since $\pi_3 = id$ on $A_3/P_3$ we have $Ker \pi_3 = P_3$. 

Suppose $\nu_{1,1}(y_2) = \sum_i y_i' \otimes y_i'' \subset A_1 \otimes A_1$. Let $x_1 \in A_1$. Denote by $\cdot$ the product in $A_{\bullet}$. Then 
$$
\nu \circ \mu(y_2\cdot x_1) = \sum_i (x_1 \cdot y_i' \otimes y_i'' + y_i' \otimes y_i''
\cdot x_1) + y_2 \otimes x_1 + x_1\otimes y_2)
$$
So a simple calculation proves the assertion.

One has
\begin{equation} \label{usy}
\pi_3(l_3(x)) = l_3(x) - \frac{1}{2} \mu(l_2(x) \otimes x) + \frac{1}{12}
\mu(l_1(x) \otimes x\otimes x)
 \end{equation}
where $l_i(x)$ is the generator of $A_n$ corresponding to the classical $n$-logarithm. For $n=2$ and $n=3$ it is given by the picture on fig.1. 
Set
$$
L_3(x):= Li_3(x) - \frac{1}{2} Li_2(x) \log x + \frac{1}{12}
Li_1(x) (\log x)^2
$$

\begin{theorem}
\begin{equation} \label{ki}
\Lambda_3(\pi_3(L,M)) = L_3(  a_3(L,M))
\end{equation}
\end{theorem}

{\bf Remark}. To get a local coincidence of two multivalued analytic functions in this formula  we should choose  appropriate cycles of integration. The theorem claims in particular that 
it is possible to do this.

{\bf Proof}. The main result of this paper implies that 
$$
\nu(\pi_3(L,M)) = \nu(l_3 \circ a_3(L,M))
$$ 
Since the differential of $\Lambda_3(L,M)$ is determined by the coproduct 
of $(L,M)$ (see the lemma below) this implies that the differentials of   both sides of (\ref{ki}) coincide. So  the difference between the left and right hand sides is a constant. Considering the additivity relation in $L$ (which has odd number of terms) we deduce that this constant is zero.

\begin{lemma}
Let $(L,M) \in A_n({\C})$ and the $A_{n-1} \otimes A_1$ component of 
  $\nu(L,M)$ is $\sum x_i \otimes y_i$. Then
$$
d \Lambda_3(L,M) = \sum d \Lambda_3(x_i) d \log y_i
$$
\end{lemma}

This is a particular case of a general fact about the differential of  the period of $n$-framed  mixed Tate motives, see the chapter "Periods" in [G9].
 A direct proof can be given by an explicit calculation using proposition 2.3.

\section{Configurations of $2n$ points in ${\Bbb P}^{n-1}$, Grassmannian $n$-logarithms  and motivic  Lie coalgebra 
of a field}

{\bf 1. Grassmannian polylogarithms}. Let me   recall the construction of the 
Grassmannian polylogarithm function ${\cal L}_n^G(h_1,...,h_{2n})$ given in [G4]. 
It is a function on configurations of {\it arbitrary} $2n$ hyperplanes in 
${\Bbb C}{\Bbb P}^{n-1}$.

Let $f_1,...,f_m$ be $m$ complex-valued functions on a manifold $X $.
 We attach to them the following  $(m-1)$-form. Let $c_{j,m}:= \frac{1}{(2j+1)!(m-2j-1)!} $.   Set
\begin {equation} \label{1}
\omega_{m-1}(f_1,..., f_m) :=
\end {equation}
$$
 \frac{1}{(2\pi i)^{m}}{\rm Alt}_m \sum_{j\geq 0} c_{j,m}\log|f_1|d\log|f_2|
\wedge ... \wedge d\log|f_{2j+1}|\wedge di\arg f_{2j+2}\wedge ... \wedge
di\arg f_{m}
$$

 Now let $f_i$ be a rational function on 
${\Bbb C} {\Bbb P}^{n-1}$ with the divisor $(f_i) = h_i - h_{2n}$, $i = 1,...,2n-1$.  
 Then
$$
{\cal L}_n^G(h_1,...,h_{2n}): = \int_{{\Bbb C} {\Bbb P}^{n-1}}\omega_{2n-2}(f_1,...,f_{2n-1})
$$

\begin{theorem} \label{gpgp}
${\cal L}_n^G(h_1,...,h_{2n})$ satisfies the following functional equations:

a) For {\it any} $2n+1$ hyperplanes $h_1,...,h_{2n+1}$ in ${\C} {\Bbb P}^{n-1}$ one has
$$
 \sum_{i=1}^{2n+1}(-1)^i{\cal L}_n^G(h_1,...,\hat h_i,...,h_{2n+1}) = 0
$$
b) For {\it any} $2n+1$ hyperplanes $p_1,...,p_{2n+1}$ in ${\C} {\Bbb P}^{n }$ one has
$$
 \sum_{i=1}^{2n+1}(-1)^i{\cal L}_n^G(p_1 \cap p_i,... ,p_{2n+1}\cap p_i) = 0
$$

c) The function ${\cal L}_n^G(h_1,...,  h_{2n })$  is skewsymmetric with respect to $h_1,...,h_{2n}$.

d)  ${\cal L}_n^G(h_1,...,  h_{2n}) =0$  if the intersection of certain $2k$
hyperplanes among $h_1,...,  h_{2n}$  has codimension $\leq k$.  
\end{theorem}

{\bf Proof}. a), b) are proved in [G4] and c) is clear.  

d). One can show that up to a (computable) rational number
$d_n$ one has
\begin{equation} \label{opopop}
\int_{{\Bbb C} {\Bbb P}^{n-1}}\omega _{2n-2}(f_1,...,f_{2n-1}) = d_n \cdot \int_{{\Bbb C} {\Bbb P}^{n-1}}\log |f_1| d \log|f_2| \wedge ... \wedge d \log|f_{2n-1}|
\end{equation}
If $h_2 \cap ... \cap h_{2k+1}$ has codimension $k$ then $d \log|f_2| \wedge 
... \wedge d \log|f_{2k+1}| =0$, so the integral is zero.

On the other hand since for $(L,M) \in A_n({\C})$ the function 
$\Lambda_n(L,M)$ is a period of the mixed Hodge structure
$ 
H^n({\C}{\Bbb P}^n \backslash L,M)
$, 
one can define  a single valued function ${\cal L}_n(L,M) \in \R$ as the $\R$-period of this mixed Hodge structure. Then ${\cal L}_n: A_n({\C})  \to \R$ is a homomorphism of groups. 

{\bf Problems}. 1. Find a complete list of functional equations for 
the function ${\cal L}_n^G$. 

2. Find the relationship between the Aomoto and Grassmannian polylogarithms. 

In the next section we will formulate these problems more precisely  and explain their importance.

{\bf 2. Conjectures}. I conjecture that one can describe {\it explicitly} a certain subgroup 
$$
{\cal R}_n^G(F)\subset \bar C_{2n}(F)$$
such that ${\cal R}_n^G({\C})$ is the group of all functional equations for the Grassmannian $n$-logarithm function ${\cal L}_n^G$. The subgroup should be given by an  explicit "universally defined"  finite list of "relations". More precisely, one should have a finite set of varieties $R_i(n)$ over ${\Z}$ and 
morphisms 
$s_i(n): {\Z}[{R}_i(n)] \to {\Z}[\overline C_{2n}({\Bbb P}^{n-1})]$ given by finite correspondences   over $\Z$ such that 
$$
{\cal R}_n^G(F) = \sum s_i(n)\Bigl({\Z}[R_i(n)]\Bigr)
$$
The properties of the subgroups ${\cal R}_n^G(F)$ are formulated in the conjecture below. Then
$$
G_n(F) = \quad \frac{{\Q}[\mbox{ configurations of any $2n$ $F$-points in ${\Bbb P}^{n-1}$} ]}{{\cal R}_n^G(F)}
$$
Then $
{\cal L}_n^G: G_n({\C})\to {\R}, \quad (l_1,...,l_{2n}) \lms {\cal L}_n^G(l_1,...,l_{2n})$. 

Let $
G_{\bullet} := \oplus^{\infty}_{n =1}G_n$. 
Recall that a structure of a Lie coalgebra on $G_{\bullet}$ is given by   homomorphisms
$$
G_n \stackrel{\delta_n}{\lra}  \oplus_{i\leq n/2} G_i \wedge G_{n-i}
$$ 
such that if $\delta:= \oplus\delta_n$ then
$$
G_{\bullet} \stackrel{\delta }{\lra}  \Lambda^2G_{\bullet}  \stackrel{\delta \otimes id - id \otimes \delta}{\lra}  \Lambda^3G_{\bullet} \lra ...
$$
is a complex.

\begin{conjecture}
a) $G_{\bullet}$ has a natural structure of a graded Lie coalgebra over ${\Q}$.

b) The category of graded finite dimensional modules over $G_{\bullet}(F)$ is equivalent to the category of mixed Tate motives over $Spec$$F$. In particular one should have
$$
H^i_{(n)}\Bigl( G_{\bullet}(F)\Bigr) = gr^{\gamma}_nK_{2n-i}(F)\otimes {\Q}
$$
\end{conjecture}
Here $H^i_{(n)}$ is the degree $n$ part of $H^i$.
 Set
\begin{equation} \label{iden}
G_1(F):= F^* \otimes {\Q},\quad G_2(F):= B_2(F) \otimes {\Q}, \quad G_3(F):= B_3(F) \otimes {\Q}
\end{equation} 
Then the first components of the cobracket $\delta $ are given by
$$
\delta_2: B_2(F) \to \Lambda^2F^*_{\Q}, \qquad \quad \delta_3: B_3(F) \to G_2(F)\otimes F^*_{\Q}
$$
where the homomorphisms $\delta_2$ and $\delta_3$ were defined in s.3.1 and 3.2. 

The $7$-term relation for the generalized cross-ratio $r_3$ is rather mysterious. Its analog 
  needed to define the group $G_4(F)$ is unknown. However the results of this paper suggest the following strategy.
Denote  by $\tilde A_n$ the free abelian group generated by the generators of the group
$ A_n$, i.e. by admissible pairs of  simplices. Let us assume that we  have defined already 
the subgroups ${\cal R}_m^G$ for $m<n$.

\begin{conjecture}There exists a homomorphism $a_n: \tilde A_n \to 
\overline C_{2n}$ such that 

a) The following diagram is commutative
$$
\begin{array}{ccc} \label {a323}
\tilde A_n&\stackrel{\nu}{\longrightarrow} &\oplus_{1 \leq i \leq n} A_i \otimes A_{n-i}   \\
&&\\
\downarrow a_n &&\downarrow a_i \wedge a_{n-i} \\
&&\\
\overline C_{2n}  &\stackrel{\delta_n}{\longrightarrow} & \oplus_{i \leq n/2}G_i(F) \wedge G_{n-i}(F)
\end{array}
$$

b) ${\cal L}_n(L;M) = c_n {\cal L}_n^G(a_n(L,M))$ for any $(L;M) \in A_n({\C})$, 
where $c_n$ is a normalization constant.
\end{conjecture}
(See also conjecture 1.42 in [G1]).

Assuming this we  {\it introduce} ${\cal R}_n^G(F)$ as   the image under the map $a_n$ of the defining relations for the group $A_n$. 
Thus $\delta_n({\cal R}_n^G)=0 $ and we  are getting a commutative diagram
 $$
\begin{array}{ccc} \label {a32323}
 A_n&\stackrel{\nu}{\longrightarrow} &\oplus_{1 \leq i \leq n} A_i \otimes A_{n-i}   \\
&&\\
\downarrow a_n &&\downarrow a_i \wedge a_{n-i} \\
&&\\
G_n(F)  &\stackrel{\delta_n}{\longrightarrow} & \oplus_{i \leq n/2}G_i(F) \wedge G_{n-i}(F)
\end{array}
$$
The results of s. 3 show how this program works  for $G_3$. 

Recall the subgroup $P_n \in A_n$ defined in the section 3.6

\begin{conjecture} \label{ih}
The map $a_n$ induces an isomorphism
$$\bar a_n: A_n/P_n \to G_n(F)
$$
\end{conjecture}
For $n=2$ this was proved in [BGSV1-2]. 

For $n=3$ one can prove that $a_3(P_3) =0$. To show that $\bar a_3$ is an isomorphism one should construct a homomorphism $L_3: B_3(F) \to A_3(F)$ splitting the map $a_3$. The homomorphism $L_3: {\Z}[F^*] \to A_n, \quad \{x\}_3 \lms L_3\{x\}_3$  is given  by the right hand side of formula (\ref{usy}). 
Then one should prove that the map $L_3: {\Z}[F^*] \to A_3/P_3$ is surjective and $L_3(R_3)=0$.

Conjecture (\ref{ih}) 
just means that the dual to the Hopf algebra $A_{\bullet} := \oplus A_n$
 is isomorphic to the universal enveloping algebra of the Lie algebra 
$G_{\bullet}^{\vee}$. (Here $A \to A^{\vee}$ is the duality between the ind and pro $\Q$-vector spaces). It seems quite remarkable that the universal enveloping algebra of $G_{\bullet}^{\vee}$ admits a completely different description (the only similar situation which comes to   mind is Lusztig's construction of $U({\cal N})$). So we get two  different descriptions of the motivic 
Lie algebra
(reflecting  the properties of  the  Aomoto and Grassmannian polylogarithms).
 It is even more interesting that there are two more  ways of thinking about 
the  same Lie algebra (!): the wonderful "cycle" 
construction of Bloch and Kriz [BK] (so far the only one which is  completely done), and  the construction reflecting  the
 properties of multiple polylogarithms ([G5]).   This definitely shows the 
richness of the subject.

All constructions of different models of the motivic Hopf algebra of the 
category of mixed Tate motives are based on the following idea. The set of appropriately 
defined equivalence classes of $n$-framed mixed Tate motives over $F$ form an abelian 
group 
${\cal A}_n$ and ${\cal A}_{\bullet}:= \oplus_{n \geq 1} {\cal A}_n$ is a 
commutative Hopf algebra. It is isomorphic to the fundamental Hopf algebra of the category of mixed Tate motives over $F$ (see [BGSV1-2], 
[BMS] for the definition of $n$-framed mixed motives). Consider a 
universal 
variation of $n$-framed mixed Tate motives  over a base $X_n$.   For any $F$-point of $x \in X_n$ we get an element 
$m_n(x) \subset {\cal A}_n$. Universality of the variation means that 
the map ${\Z}[X_n(F)] \to {\cal A}_n$ is surjective. The kernel of 
this map is supposed to be described explicitly. So 
$$
\oplus_{n >0} \frac{{\Z}[X_n(F)]}{{\rm Ker } \quad m_n}
$$ 
should have a natural structure of a Hopf algebra, and
one needs to determine it. (The variation of mixed Tate motives for the ``cycle'' Hopf algebra of [BK] can be found in the last section of [G4]).  
So it is quite interesting that we could get a 
construction of a co-Lie algebra directly, without constructing first its universal enveloping algebra. 

{\bf 3.  A canonical map ${\Z}[F^*] \to A_n(F)$}. Let $A$ be a commutative graded  Hopf algebra with a product $\mu$ and coproduct $\nu$, $A_+$  the kernel of the augmentation homomorphism and  $$
\tilde \nu:= \nu - (id \otimes 1 + 1 \otimes id): A \lms A_+^{\otimes 2}
$$
 the restricted coproduct. We define a map of graded vector spaces
$\tilde \nu_{[k]}: A \lra A_+^{\otimes k}$ as a composition 
$$
A \stackrel{\tilde \nu}{\lra} A_+\otimes A_+ \stackrel{\tilde \nu \otimes id}{\lra} A_+\otimes A_+ \otimes A_+ \stackrel{\tilde \nu \otimes id}{\lra} ... \stackrel{\tilde \nu \otimes id}{\lra} A_+^{\otimes k}
$$
  Let $\mu_k:A^{\otimes k} \to A$ be the product map. Set
$$
\pi:= \sum_{k\geq 1}\frac{(-1)^{k-1}}{k}\mu_k \circ \tilde \nu_{[k]} = Id - \frac{1}{2} \mu \circ \tilde \nu + ... 
$$
(The map $\pi_3$ from the section 3.6 is a particular case of this)
\begin{proposition} \label{hhoo}
$\pi^2 = \pi$ and $Ker \pi =P$.
\end{proposition}

{\bf Proof}. Since $\pi =id $ (mod $P$) one needs to show only that $P \subset Ker \pi$. 

Now return to the Hopf algebra $A_{\bullet}(F)$. Let $B_k$ be the Bernoulli numbers. Recall that $l_n(x)$ is an element of $A_n(F)$ corresponding to the classical $n$-logarithm. 

\begin{proposition} \label{o}
$$
\pi (l_n(x)) = \sum_{k \geq 0}\frac{B_k}{k!}l_{n-k}(x)\cdot x^k
$$ 
\end{proposition}
Here $x^m:= x \cdot ... \cdot x$ ($m$ times) is the product in $A_{\bullet}(F)$ of the element of $A_1$ corresponding to $x$ under the canonical isomorphism $A_1(F) \to F^*$. Notice that  this formula  coincides with  the formula for the function $\Lambda_n(z)$ in s. 4.1 of [BD]. However we get it in a quite different way. 

{\bf Proof}. 
Let
$$
l(x,t):= \sum_{k \geq 1} l_k(x)t^{k} 
$$
\begin{lemma} \label{stand}
$\tilde \nu: l(x,t) \lms l(x,t) \otimes (e^{x\cdot t} -1)$
\end{lemma}

This is the  generating function for the standard formula 
$$
\tilde \nu (l_k(x)) = \quad \sum_{1 \leq i\leq k-1}l_{k-i}(x) \otimes \frac{x^i}{i!}
$$
for the coproduct of the classical polylogarithm. 
Therefore 
$$
\mu \circ \nu_{[k]}: \quad l(x,t) \lms l(x,t) \cdot (e^{x\cdot t} -1)^{k-1}
$$ 
So
$$
\pi: l(x,t) \lms \quad - \sum_{k \geq 1} \frac{(-1)^k}{k}l(x,t) 
\cdot(e^{x\cdot t} -1)^{k-1} = 
$$
$$
=l(x,t) \sum_{k \geq 1} \frac{(1-e^{x\cdot t})^k}{k}\frac{1}{e^{x\cdot t} -1}\quad = \quad l(x,t) \cdot\frac{xt}{e^{x\cdot t} -1} \quad = \quad \sum_{n \geq 1}\Bigl(\sum_{k \geq 0}\frac{B_k}{k!}l_{n-k}(x)\cdot x^k \Bigr) t^n
$$
The proposition is proved.

The kernel of the map $$
 L_n: {\Z}[F^*] \lms A_n(F), \quad \{x\} \lms \pi (l_n(x))
$$ 
should coincide with the subgroup of all functional equations for the $n$-logarithm.

 \section  {Motivic  structure of the Grassmannian tetralogarithm and Lie coalgebra $G(F)_{\leq 4}$}

{\bf 1. The group $\tilde G_{2n}(F)$}.
 Let $\tilde G_{2n}(F)$ be the free abelian group generated by 2n-tuples of
points $(l_1,...,l_{2n})$ in generic position in ${{\Bbb P}}^{n-1}(F)$ subject to the following relations:

1){\it Projective invariance}: $(l_1,...,l_{2n}) = (gl_1,...,gl_{2n})$ for any $g \in PGL_{n}(F)$.

2) {\it Skew symmetry}. $(l_1,...,l_{2n} ) = (-1)^{|\sigma|} (l_{\sigma
(1)},...,l_{\sigma (2n)})$
for any   $\sigma \in S_{2n}$.

3) {\it $(2n+1)$ -term relation}: for any $2n+1$ points in generic
position $(l_0,...,l_{2n})$ in $P^{n-1}(F)$ one has
$ 
\sum_{i=0}^{2n}(-1)^{i}(l_0,...\hat
l_i,...,,l_{2n}) =0
$ 

4){\it dual $(2n+1)$ -term relation}: for any $2n+1$ points in generic
position $(l_0,...,l_{2n})$ in ${{\Bbb P}}^{n}(F)$ one has
$ 
\sum_{i=0}^{2n}(-1)^{i}(l_i|l_0,...\hat
l_i,...,,l_{2n}) =0
$ 

{\bf Remark}. These  relations   reflect  the properties of the Grassmannian $n$-logarithm listed in theorem 4.1. The group $G_{2n}(F)$ is supposed to be a quotient of the group $\tilde G_{2n}(F)$. It is a nontrivial quotient already for $n=3$. 

{\bf 2. The main result}. After identification (\ref{iden})
  the degree $4$ part of the cochain complex of the Lie coalgebra $G_{\bullet}(F)$ should look  as follows: 
$$
G_4(F)  \quad \stackrel{\delta}{\lra} \quad  B_3(F) \otimes F^{\ast}
\quad \oplus \quad   B_2(F)\wedge  B_2(F)  \stackrel{\delta}{\lra} B_2(F) \otimes \Lambda^2F^{\ast} \stackrel{\delta}{\lra} \Lambda^4F^{\ast} 
$$
Here  
$$
\delta: \{x\}_3 \otimes y \lms \{x\}_2 \otimes x \wedge y, \qquad  
$$
$$
\delta: \{x\}_2 \wedge \{y\}_2 \lms\{y\}_2 \otimes  (1-x) \wedge x - \{x\}_2 \otimes  (1-y) \wedge y
$$

 It remains to define   a homomorphism 
$$
\tilde G_4(F) \quad \stackrel{ \delta }{\longrightarrow} \quad B_3(F) \otimes F^{\ast}
\quad \oplus  \quad B_2(F)\wedge  B_2(F)   
$$
  Below we will   construct it  as a composition
$$
\tilde G_4(F) \quad \stackrel{ \delta }{\longrightarrow}\quad  A^0_3(F) \otimes F^{\ast} \quad \oplus \quad F^{\ast}  \otimes A^0_3(F)  
\quad \oplus \quad B_2(F) \wedge  B_2(F)  
$$
$$
\stackrel{(a_3 \otimes id - id \otimes a_3
, Id)}{\longrightarrow}B_3(F) \otimes F^{\ast} \quad \oplus \quad  B_2(F)\wedge  B_2(F)   
$$

Namely, $ \delta = ( \delta_{3,1},   \delta_{1,3},  \delta_{2,2})$ where
$$
 \delta_{3,1}(l_1,...,l_8) : = {\rm Alt}_8
\Bigl((l_1,l_2,l_3,l_4;l_5,l_6,l_7,l_8)  \otimes \Delta (l_5,l_6,l_7,l_8)
\Bigr) \in A^0_3(F) \otimes F^*
$$
$$
 \delta_{1,3}(l_1,...,l_8) : = -{\rm Alt}_8
\Bigl(\Delta (l_1,l_2,l_3, l_4) \otimes (l_1,l_2,l_3,l_4;l_5,l_6,l_7,l_8) 
\Bigr) \in F^* \otimes A^0_3(F)
$$
$$
 \delta_{2,2}(l_1,...,l_8) : =  \frac{288}{7} \cdot {\rm Alt}_8
\Bigl( (l_1,l_2|l_3,l_4,l_5,l_6)_2\wedge (l_3, l_4|l_1,l_2,l_5,l_7)_2\Bigr) \in \Lambda^2 B_2(F)
$$

\begin{theorem} \label{mmmar} a) The homomorphism $ \delta$ is well defined, i.e. sends the relations 1) - 4) to zero.

b) The composition
$$
\tilde G_4(F) \stackrel{ \delta }{\longrightarrow} B_3(F) \otimes F^{\ast}
\quad \oplus \quad  B_2(F)\wedge  B_2(F)  \stackrel{\delta}{\lra} B_2(F) \otimes \Lambda^2F^{\ast} 
$$
is zero.
\end{theorem}
  
The proof is postponed to Section 5.6.

{\bf Remark}. Taking, as usual ([G1]), the "connected component" of zero in ${\rm Ker} \delta $ we should get the set of defining relations for the group $G_4(F)$. However an explicit construction of them is not known yet.

 {\bf 3. Applications to the Borel regulator map}. Recall the rank filtration quotient $K^{[3]}_7(F)\otimes {\Q} $ of $K_7(F)\otimes {\Q} $,  which is expected to be isomorphic to $gr^{\gamma}_4K_7(F)\otimes {\Q}$, see [G1]. 

It turn  out that using the  results of section 3.3, namely the definition of the map $a_3$ given there  and theorem 3.3, we can factorize the natural 
projection of the map $ \delta$ to $B_3(F) \otimes  F^* \oplus \Lambda^2 B_2(F)$, i.e. the map
$$
  C_8(V_4) \lra B_2(F) \otimes \Lambda^2F^*
$$ 
 as a composition
$$
C_8(V_4) \stackrel{\partial}{\lra} C_7(V_4) \stackrel{f_7(4) }{\lra} B_3(F) \otimes  F^* \oplus \Lambda^2 B_2(F)$$
where the  map   
$$
f_7(4): C_7(V_4)  \lra B_3(F) \otimes F^* \oplus \Lambda^2B_2(F)
$$
is defined as follows. 

We will use  shorthands like : 
$$
( 1,  2, 3, 4) \quad \mbox{for} \quad  \Delta(l_1,l_2,l_3,l_4)\quad \mbox{and so on}
$$ 
Write $f_7(4) = f_{2,2} + f'_{3,1}+ f''_{3,1}$, where 
$$
f_{2,2}: (l_1,...,l_7) \lms   \frac{288}{7} \cdot {\rm Alt}_7
\Bigl( (l_1,l_2|l_3,l_4,l_5,l_6)_2\wedge (l_3, l_4|l_1,l_2,l_5,l_7)_2\Bigr) \in \Lambda^2 B_2(F)
$$
$$
f'_{3,1}: (l_1,...,l_7) \lms \quad \frac{-32}{6} \cdot {\rm Alt}_7
\Bigl( (l_1| l_2,l_3,l_4,l_5,l_6,l_7)_3 \otimes \Delta(l_1,l_2,l_3,l_4)\Bigr) \in  B_3(F) \otimes F^*
$$
$$
f''_{3,1}: (l_1,...,l_7) \lms \quad 96 \cdot {\rm Alt}_7
\Bigl( \{r (l_1,l_5| l_2,l_6,l_3,l_7)\}_3  \otimes \Delta(l_1,l_2,l_3,l_4)\Bigr)\in  B_3(F) \otimes F^* 
$$
Notice that if we divide all the coefficients by $32$ the coefficients will be smaller:  $\frac{9}{7}, \frac{-1}{6}, 3 $. We will use a notation $f_{3,1}$ for $f_{3,1}' + f_{3,1}^{''}$. Finally, $(l_1| l_2,l_3,l_4,l_5,l_6,l_7)_3 \in B_3(F)$ 
is obtained as follows: the configuration $(l_1| l_2,l_3,l_4,l_5,l_6,l_7)$ of $6$ points in $P^2$  provides a generator of $G_3(F)$, which is then mapped to $B_3(F)$ using the generalized cross-ratio $r_3$. 

To get the formula for $f''_{3,1}$ we have used the definition of $a_3''$ and the intermediate formula
$$
f''_{3,1} \circ \partial: (l_1,...,l_8) \lms \quad \frac{-32}{3}  \cdot {\rm Alt}_8
\Bigl( \mu_3(l_1,  l_5| l_2, l_3,l_4; l_6,l_7, l_8)  \otimes \Delta(l_1,l_2,l_3,l_4)\Bigr)
$$

Therefore there is the following commutative diagram:
$$
\begin{array}{ccccccc}
C_9(V_4) &\stackrel{\partial }{\lra}&C_8(V_4) &\stackrel{\partial }{\lra}&C_7(V_4)&& \\ 
&&&&&&\\
\downarrow &&\downarrow f_8(4)&&\downarrow f_7(4)&&\\
&&&&&&\\
0& \lra & \tilde G_4(F)& \stackrel{\delta}{\lra}& B_3(F) \otimes F^* \oplus \Lambda^2B_2(F)& \lra &B_2(F) \otimes \Lambda^2F^*\\
\end{array}
$$
Moreover, it is easy to see that the composition 
$$
C_9(V_5)  \stackrel{\partial'}{\lra}  C_8(V_4)  \lra B_3(F) \otimes F^* \oplus \Lambda^2B_2(F)
$$
 is zero. Here $(l_1,...,l_9) \lms \sum(-1)^i (l_i|l_1,..., \widehat
l_i,...,l_9)$. 

Therefore we have constructed a morphism from the appropriate part of the weight $4$ bigrassmannian complex ([G6]) to the bottom line in the diagram above. (The full homomorphism from the weight $4$ bigrassmannian complex 
to the weight 4 motivic complex will be treated in [G8]). 
Applying the general  technique developed in [G1-2]  and  [G6] we get
  part a) of the following theorem. 

\begin{theorem} \label{5.2th} a) There exists a canonical  map 
$$
K^{[3]}_7(F)\otimes {\Q} \quad \lra \quad {\rm Ker} \Bigl(\tilde G_4(F)  \stackrel{\delta}{\lra} B_3(F) \otimes F^* \oplus \Lambda^2B_2(F)\Bigr)_{\Q}
$$

b) In the case $F = {\C}$ 
the composition 
$$
K^{[3]}_7({\C} ) \lra \tilde G_4({\C})  \stackrel{{\cal L}^G_4}{\lra} \R
$$
coincides with a nonzero rational multiple of the Borel regulator map. 
\end{theorem}

Namely, the diagram above provides a morphism from a piece of the
weight four bigrassmannian
complex to the complex $$
0 \lra \tilde G_4(F)_{\Q}  \stackrel{\delta}{\lra}
(B_3(F) \otimes F^* \oplus \Lambda^2B_2(F))_{\Q}
$$
 Combining this
with the canonical maps from the homology of $GL(F)$ to the weight
four complex of
affine flags ([G6]) followed by the canonical map from the complex of
affine flags to the 
weight four bigrassmannian
complex we get a canonical map 
$$
H_7(GL(F), {\Q}) \quad \lra \quad {\rm Ker} \Bigl(\tilde G_4(F)
\stackrel{\delta}{\lra} B_3(F) \otimes F^* \oplus
\Lambda^2B_2(F)\Bigr)_{\Q}
$$ Using the arguments given in the 
papers cited above we get   part a).

To get  part b)  we use the
 computation of the Borel regulator map via the Grassmannian
$n$-logarithm ${\cal L}_n^G$  (see  [G4] and [G7]). 

\begin{conjecture} There exists a map 
$$
{\rm Ker} \Bigl(\tilde G_4(F)  \stackrel{\delta}{\lra} B_3(F) \otimes F^* \oplus \Lambda^2B_2(F)\Bigr)_{\Q} \quad \lra \quad K^{[3]}_7(F)\otimes {\Q} 
$$
 which in the case $F = {\C}$ commutes with the Borel regulator map. 
\end{conjecture}

{\bf 4. Towards Zagier's conjecture on $\zeta_F(4)$}. Let us    construct  a homomorphism 
\begin{equation} \label{3.62}
\overline f_7(4): C_7(4) \longrightarrow  B_{3}(F)\otimes F^{\ast}
\end{equation}
providing a definition of a homomorphism  $\overline \delta: G_{4}(F) \lra B_3(F) \otimes F^*$ whose composition with the natural map $B_3(F) \otimes F^* \lra B_2(F) \otimes \Lambda^2F^*$ is zero. So we will get a commutative diagram 
\begin{center} 
\begin{picture}(150,70) 
\put(-40,50){$C_8(4)$}  
\put(-45,0){$ G_{4}(F)$} 
\put(70,50){$ C_7(4)$} 
\put(70,0){$B_{3}\otimes F^{\ast}$} 
\put(-20,45){\vector(0,-1){30}} 
\put(-44,25){$f_8(4)$}  
\put(5,50){\vector(1,0){55}}  
\put(80,45){\vector(0,-1){30}} 
\put(85,25){$\overline f_7(4)$} 
\put(5,5){\vector(1,0){55}} 
\put(30,10){$\overline \delta$} 
\put(30,60){$\partial $}  
\put(250,35){} 
\end{picture} 
\end{center}

This is done in 2 steps. First we use the formulas for the map $  f_7(4)$ given in s. 5.3 as a  definition of   a
homomorphism 
\begin{equation} \label{**}
C_7(V_4) \stackrel{f_{3,1} \oplus f_{2,2} }{\longrightarrow} B_{3}\otimes
F^{\ast} \oplus \Lambda^2 {\Bbb Z}[P^1_F]
\end{equation} 

 Then we construct  a homomorphism  
$$
g: \Lambda^2 {\Bbb Z}[P^1_F] \longrightarrow B_{3}\otimes F^{\ast} 
$$
making the diagram
$$
\begin{array}{ccc}
 \Lambda^2 {\Bbb Z}[P^1_F]&& \\
& {\searrow} \tilde \delta & \\
{\downarrow} g&&B_{2}\otimes \Lambda^{2}F^{\ast}\\
& {\nearrow}\delta & \\
B_{3}\otimes F^{\ast}&&  
\end{array}
$$
commutative, and so providing a commutative diagram 
\begin{center} 
\begin{picture}(150,70) 
\put(-96,50){$B_{3}\otimes F^{\ast} \oplus \Lambda^2 {\Bbb Z}[P^1_F] $}

\put(-65,0){$B_{3}\otimes F^{\ast}$} 
\put(70,50){$B_{2}\otimes \Lambda^{2}F^{\ast}$} 
\put(70,0){$B_{2}\otimes \Lambda^{2}F^{\ast}$} 
\put(-45,45){\vector(0,-1){30}} 
\put(-75,25){$id \oplus g$}
\put(5,50){\vector(1,0){55}} 
\put(80,45){\vector(0,-1){30}} 
\put(85,25){$id$}

\put(5,5){\vector(1,0){55}} 
\put(30,10){$\delta$}

\put(30,60){$\delta + \tilde \delta$} 
\put(250,35){$\ast $} 
\end{picture} 
\end{center}  

Here 
$$
\tilde \delta \Bigl(   \{x\} \wedge \{y\} \Bigl) = 
  \{y\}_2 \otimes (1-x) \wedge x - 
\{x\}_2 \otimes (1-y) \wedge y
$$

To define $g$ recall that the group $S_3$ acts naturally on $P^1 \backslash \{0,1,
\infty \}$.  The orbit of a point $x$
is
$\quad x, \quad x^{-1},\quad 1-x,\quad (1-x)^{-1},\quad 1-x^{-1} , \quad
(1-x^{-1})^{-1} $.
Set 
\begin{equation} \label{3.9}
g: \quad \{x\}_2 \wedge \{y\}_2 \longrightarrow -\frac{1}{12}\sum_{\sigma_1,\sigma_2 \in S_3}(-1)^{\vert\sigma_1\vert + \vert\sigma_2\vert}\{\frac{\sigma_1(x)}{\sigma_2(y)}\}_3 \otimes \frac{1- \sigma_1(x)
}{1-\sigma_2(y)}
\end{equation}
Notice that $\{x^{-1}\}_3 = \{x\}_3$ and so the right-hand side is
skew-symmetric with respect to transposition of $x$ and $y$.

We define the desired homomorphism
 (\ref{3.62}) as the composition of  morphismes  (\ref{**}) and (${\rm id} \oplus g $).  
 
 {\bf Remark}.  $g$ does  not factorizes through a  homomorphism of $\Lambda^2
B_2(F)$. In fact we proved in section 4 of [G2] that there is no  ``natural'' (i.e. given by
formulas) non-zero
homomorphism  $\Lambda^2 B_2(F) \longrightarrow B_3(F) \otimes
F^{\ast}$.

\begin {proposition} 
  One has $ \delta \circ    g  = \tilde \delta$. Therefore the diagram $(\ast  )$ is commutative. 
\end {proposition}

{\bf Proof}. Computing the coboundary $\delta$ of expression (\ref{3.9})
we get 

\begin{equation} \label{3.9x}
  \frac{1}{12}\sum_{\sigma_1,\sigma_2 \in S_3}(-1)^{\vert\sigma_1\vert + \vert\sigma_2\vert}\{\frac{\sigma_1(x)}{\sigma_2(y)}\}_2 \otimes \frac{\sigma_1(x)}{\sigma_2(y)}\wedge \frac{1- \sigma_1(x)
}{1-\sigma_2(y)} \in B_2(F) \otimes \Lambda^2F^{\ast}
\end{equation}

$\sigma_1(x) \wedge (1-\sigma_1(x)) \in \Lambda^2F^{\ast}$ is
skewsymmetric with respect to  $S_3$. So $(1-x) \wedge x
$ appears in (\ref{3.9x}) with factor
\begin{equation} \label{3.9a}
- \frac{1}{12} \sum_{\sigma_1, \sigma_2 \in S_3}(-1)^{\vert\sigma_2\vert} 
\{\frac{\sigma_1(x)}{\sigma_2(y)}\}_2 
\end{equation}

Modulo 2-torsion one can rewrite the basic
relation in $B_2(F)$ as follows
\begin{equation} \label{3.9bs}
\{x\}_2 - \{y\}_2 = \{\frac{x}{y}\} - \{\frac{1-x}{1-y}\} +
\{\frac{1-x^{-1}}{1-y^{-1}}\} 
\end{equation}

Averaging with respect to $x$ over the group $S_3$ we get 
\begin{equation} \label{3.9w}
-6 \{y\}_2 = \sum_{\sigma_1 \in S_3}
\{\frac{\sigma_1(x)}{y}\}_2 - \{\frac{\sigma_1(x)}{1-y}\}_2 + 
\{\frac{\sigma_1(x)}{1-y^{-1}}\}_2
\end{equation}
 So modulo 6-torsion (\ref{3.9a}) coincides with $ \frac{1}{2}( \{y\}_2 -
\{\frac{y}{y-1}\}_2) =  \{y\}_2$. The proposition is proved. 

   {\it A definition of the  homomorphism
$\overline \delta: G_4(F) \lra B_3(F) \otimes
F^{\ast}$}.  Take any $8$ points $(l_1,...,l_8)$ 
in generic position in ${{\Bbb P}}^3(F)$. Lift them to $8$ vectors $(\tilde l_1,...,\tilde l_8)$ in the $4$-dimensional vector space $V_4$ and then apply to them the homomorphism $\overline f_7(4) \circ \partial$. We claim that the result does not depend on the choise of vectors $\tilde l_i$ projecting to the points $l_i$. Indeed, only $f_{3,1}$ component of the map $\overline f_7(4)$ 
may depend on this choise, and it is easy to check it does not. 

Recall the group ${\cal B}_4(F)$ defined in  [G1], [G2].

\begin{conjecture} \label{uucc}
There exists a canonical homomorphism of groups ${\Bbb L}_4: G_4(F) \lra {\cal B}_4(F)$ making the following diagram commutative:
$$
\begin{array}{ccc}
 G_4(F) &\stackrel{\overline \delta}{\lra} & B_3(F) \otimes
F^{\ast}\\
&&\\
{\Bbb L}_4 \downarrow && \downarrow =\\
&&\\
{\cal B}_4(F)&\stackrel{ \delta}{\lra} & B_3(F) \otimes
F^{\ast}  
\end{array}
$$
\end{conjecture}

It follows from theorem \ref{5.2th} that this conjecture implies Zagier's conjecture for $\zeta_F(4)$ for any number filed $F$.

{\bf 5. A  motivic  construction of the Grassmannian tetralogarithm}. 
Theorem \ref{mmmar}
 allows us to construct a tetralogarithm function $\tilde {\cal L}_4^{G}$ 
on configurations of $8$ points in ${\C} {{\Bbb P}}^3$ providing    a homomorphism 
$ 
  \tilde G_4({\C}) \lra \R
$.
Namely, let $X$ be a variety over ${\C}$ and $F:= {\C}(X)$. Then there is a homomorphism of complexes
$$
\begin{array}{ccccc}
 B_3(F) \otimes F^{\ast}
\quad \oplus \quad  B_2(F)\wedge  B_2(F)& \stackrel{\delta }{\longrightarrow} &B_2(F) \otimes \Lambda^2F^*& \stackrel{\delta }{\longrightarrow}& \Lambda^4F^*\\
&&&&\\
\downarrow r_4(2)&&\downarrow r_4(3)&&\downarrow r_4(4)\\
&&&&\\
S^1(Spec F)& \stackrel{d}{\longrightarrow}&S^2(Spec F)&\stackrel{d}{\longrightarrow}& S^3(Spec F)
\end{array}
$$
where $S^{\bullet}(Spec F)$ is the de Rham complex of smooth forms at the generic point of $X$ over ${\C}$, given by the following formulas. Set
$$
\widehat{{\cal L}_{n}} (z) = \left\{ \begin{array}{ll}
{\cal L}_{n}(z) & n :\ {\rm odd} \\

i {\cal L}_{n}(z) & n: \ {\rm even} \end{array} \right.
$$
and 
$$
\alpha(g_1,g_2):= -\log|g_1| d \log |g_2| +  \log |g_2| d \log|g_1| 
$$
We define homomorphisms $r_4(\bullet)$ by the formulas
$$
\begin{array}{ccc}
 & &r_{4}(2) : \{ f \}_{3} \otimes g \mapsto  \widehat{{\cal L}}_{3} (f ) d i \arg g    -\frac{1}{3} \widehat{{\cal L}}_{2} (f )\log \vert g \vert \cdot d \log \vert f\vert 
\end{array}
$$
$$
\begin{array}{ccc}
 & &r_{4}(2) : \{ f \}_{2} \wedge \{ g \}_{2}  \mapsto   \frac{1}{3}\cdot \Bigl(   \widehat{{\cal L}}_{2} (g )   \cdot \alpha(1-f,f) - \widehat{{\cal L}}_{2} (f )   \cdot \alpha(1-g,g) \Bigr)
\end{array}
$$
$$
\begin{array}{ccc}
& & r_{4}(3) : \{ f \}_{2} \otimes g_1 \wedge g_2 \mapsto   \widehat{{\cal L}}_{2} (f ) d i\arg g_1 \wedge d i\arg g_2   -
\frac{1}{3} \alpha(1-f,f) \wedge \\
& &\Bigl(\log \vert g_1 \vert  d \arg \vert g_2\vert  -
\log \vert g_2 \vert  d \arg \vert g_1\vert\Bigr) \quad 
+ \frac{1}{3} \widehat{{\cal L}}_{2} (f ) d \log \vert g_1 \vert \wedge d \log \vert g_2 \vert\\
\end{array}
$$
 $$
r_4(4): f_1 \wedge ... \wedge f_4  \lms  \omega_3(f_1 \wedge ... \wedge f_4)
$$ 
A direct computation shows that we get a morphism of complexes.

It follows from the theorem that the composition 
$ 
  r_4(2) \circ \delta(l_1,...,l_8)  
$ 
 is  a closed 1-form on the space of generic configurations of $8$ points in ${{\Bbb P}}^3$. It turns out that integrating it we get a (single-valued) function $\tilde {\cal L}_4^{G}$:

\begin{proposition} \label{svf}
The  integral $\int_{\gamma} r_4(2) \circ \delta(l_1,...,l_8)$ is a single-valued function defined on the space of configurations of $8$ points in ${{\Bbb P}}^3$ in generic position. 
\end{proposition}

Here we integrate along a path $\gamma$ from a given reference point  to a variable point in the configuration space. The constant of integration is normalized by the condition that the function tends to zero when the configuration degenerates.  

{\bf Proof}. 
The fundamental group of the configuration space is generated by loops around divisors $\Delta(l_{i_1}, ..., l_{i_4})=0$. It is clear from the formula for $r_4(2)(\{f\}_2 \wedge \{g\}_2)$ that the $B_2 \wedge B_2$ part of the 1-form $r_4(2) \circ \delta(l_1,...,l_8)$ has trivial monodromy. Further, the monodromy around $\Delta(l_{1}, ..., l_{4})=0$ is equal to the limit value of $2 \cdot 2 \pi i \widehat {\cal L}_3 \circ a_3((l_1, l_2, l_3, l_4; l_5, l_6, l_7, l_8))$  at the divisor $\Delta(l_{1}, ..., l_{4})=0$, which is zero. The proposition is proved. 

We have constructed two functions, $\tilde {\cal L}_4^{G}$ and ${\cal L}_4^{G}$,  on  configurations of 8 points in ${\Bbb P}^3$ satisfying the properties listed in theorem 4.1.

\begin{question} \label{mmmar1}
Does the function $\tilde {\cal L}_4^{G}$  coincides with a  multiple of the Grassmannian tetralogarithm ${\cal L}_4^{G}$ constructed in section 3. More precisely,
$\tilde {\cal L}_4^{G} = \lambda(2 \pi )^3 \cdot {\cal L}_4^{G}$ where $\lambda \in {\Q}^*$?
\end{question}

I am completely sure these two functions essentially coinside. 
           
{\it A version of  conjecture \ref{uucc}}. The composition
$$
G_4({\C}) \stackrel{\overline \delta}{\lra} B_3({\C}) \otimes {\C}^* \stackrel{r_4(2)}{\lra} \quad \mbox{$1$-forms}
$$
provide {\it another}  $1$-form on the configuration space of $8$ points. 
It is closed by the main result of this section. Therefore we can integrate it and get a function, which turns out to be single-valued (the same argument as above),  denoted $\overline {\cal L}^G_4$. 

Conjecture  \ref{uucc} implies that this function is expressable via the classical $4$-logarithm.

{\bf 6. Proof of theorem \ref{mmmar}}.  It consists of two different parts. 
We will first compute the image in $B_2 \otimes \Lambda^2F^*$ of the $\delta_{1,3} + \delta_{3,1}$-component of $\delta$, then do the same for the $\delta_{2,2}$-component, and will see that they differ by a sign.

{\it Part 1}. During the proof of this theorem we will use  shorthands like : 
$$
( 1,  2| 3, 4, 5, 6)_2 \quad \mbox{for} \quad \{r(l_1,l_2|l_3,l_4,l_5,l_6)\}_2
$$ 

\begin{lemma}
The following composition  
$$
 \tilde G_4 \stackrel{\delta_{3,1}}{\lra} A_3^0 \otimes F^* \stackrel{\nu_{1,2}\otimes {\rm id}}{\lra} F^*\otimes A_2 \otimes F^* \stackrel{p}{\lra} B_2 \otimes \Lambda^2F^*,
$$
where $p: x_1 \otimes y_2 \otimes z_1 \lms -a_2(y_2)\otimes x_1 \wedge z_1$, is equal to zero. 
\end{lemma}

{\bf Proof}. Using proposition 2.3 
$$
p \circ (\nu_{1,2}\otimes {\rm id}) \circ \delta_{3,1}: \quad (l_1,...,l_8) \lms
$$
$$
-16 \cdot {\rm Alt}_8   \Bigl( a_2 (5|2,3,4;  6,7,8) 
       \otimes    (5,2,3,4 )  \wedge ( 5,6,7,8)\Bigr) = 0
$$
Indeed, the expression we alternate is symmetric with respect to the odd involution exchanging $(2,3,4)$ with $(6,7,8)$. Here we used the formula 
$a_2(L,M) = - a_2(M,L) $ ([BGSV]).

Let us compute the composition
\begin{equation} \label{kk1}
G_4 \stackrel{\delta_{3,1}}{\lra} A_3^0 \otimes F^* \stackrel{\nu_{2,1}\otimes {\rm id}}{\lra} A_2 \otimes F^*\otimes F^* \stackrel{  }{\lra} B_2 \otimes \Lambda^2F^*
\end{equation}

\begin{lemma}
The composition (\ref{kk1}) is given by
$$
(l_1,...,l_8) \lms - 144 \cdot {\rm Alt}_8  \Bigl(  (3,4|1,2, 6,7)_2
       \otimes    ( 1,2,3,4)  \wedge ( 1,2,3,5)\Bigr)
$$
\end{lemma}

{\bf Proof}. Using proposition 2.3 we get
$$
\nu_{2,1}\otimes {\rm id} \circ \delta_{3,1}: \quad (l_1,...,l_8) \lms
$$
$$
16 \cdot {\rm Alt}_8  \Bigl(  (1|2,3,4;  6,7,8) 
       \otimes    ( 1, 6,7,8)  \wedge ( 5,6,7,8)\Bigr)
$$
Applying formula (6) for $a_2$ and using $a_2(L,M) = -a_2(M,L) $ we get 
$$
- 144 \cdot {\rm Alt}_8  \Bigl(  (1,6|7,8,3,4)_2 
       \otimes    ( 1,6,7,8)  \wedge (5,6,7,8)\Bigr)
$$
Using the even permutation $(1,2,3,4,5,6,7,8) \lms (4,8,6,7,5,3,1,2)$  
we get the lemma.

The computations for $\delta_{1,3}$ are  completely similar and can be formally deduced from the computations for $\delta_{3,1}$ using the following fact. Consider the composition
$$
A_3^0 \quad \stackrel{\nu_3}{\lra} \quad A_2 \otimes  A_1 \quad \oplus \quad A_1 \otimes  A_2 \quad 
\stackrel{a_2 \wedge a_1}{\lra }\quad   B_2 \otimes F^*
$$
where  $a_2 \wedge a_1$ was defined in s. 3.3. Then $a_2 \wedge a_1\circ \nu_3(L,M) = a_2 \wedge a_1\circ \nu_3(M,L)$.

{\it Part 2}.
We will use a lot the five term relation
$\sum(-1)^i (l_i|l_1,...,\widehat l_i, ... ,l_5)_2 =0$.  

Since
$$
\delta (1,2|3,4,5,6)_2 = \frac{1}{2}{\rm Alt}_{\{3,4,5,6\}} ((1,2,3,4)  \wedge (1,2,3,5))
$$
 and 
$$
{\rm Alt}_8  \Bigl(  (1,2|3,4,5,6)_2
           \wedge (3,4|1,2,5,7)_2\Bigr) = 
$$
$$
-{\rm Alt}_8  \Bigl(
          (3,4|1,2,5,7)_2 \wedge ( (1,2|3,4,5,6)_2\Bigr)
$$
one has
$$
\delta \circ {\rm Alt}_8  \Bigl(( (1,2|3,4,5,6)_2
           \wedge (3,4|1,2,5,7)_2\Bigr) =
$$
$$
=  {\rm Alt}_8\Bigl((3,4|1,2,5,7)_2 \otimes  {\rm Alt}_{\left\{3,4,5,6 \right \}}
 [(1,2,3,4)  \wedge (1,2,3,5) ]\Bigr)  = 
$$
$$
 {\rm Alt}_8\Bigl((3,4|1,2,5,7)_2 \otimes 
 [ -2\cdot (1,2,3,4)  \wedge (1,2,3,6)  + 2\cdot (1,2,3,5)  \wedge (1,2,3,6) - 
$$
$$
 2\cdot (1,2,3,5)  \wedge (1,2,5,6)  + 2\cdot (1,2,3,6)  \wedge (1,2,5,6)  -  (1,2,3,6)  \wedge (1,2,4,6) ]\Bigr)  
$$
In the last step we have used the following simple observation: the terms in the $\Lambda^2F^*$-factor where   $6$ is absent vanish by skew-symmetry,  since $6$ and $8$  are also absent in  $(3,4|1,2,5,7)_2$. 

 Using the skew-symmetry we rewrite the last formula as follows:  
\begin{equation} \label{001}
 {\rm Alt}_8 \Bigl(\Bigl(  2 \cdot (3,4|1,2,6,7)_2
+ 2\cdot (3,6|1,2,5,7)_2    + 2\cdot (4,6|1,2,3,7)_2 
\end{equation}
$$
  + 2\cdot (4,6|1,2,5,7)_2 + (5,4|1,2,6,7)_2 \Bigr)
\otimes   (1,2,3,4)  \wedge (1,2,3,5)  \Bigr)
$$
\begin{lemma} \label{00}
The  expression (\ref{001}) is equal to
$$
 7 \cdot {\rm Alt}_8\Bigl(  (3,4|1,2,6,7)_2  
\otimes   (1,2,3,4)  \wedge (1,2,3,5)  \Bigr)
$$
\end{lemma}

{\bf Proof}. Computing 
 \begin{equation} \label{ram1}
{\rm Alt}_8\Bigl(   2\cdot (3,6|1,2,5,7)_2   
\otimes   (1,2,3,4)  \wedge (1,2,3,5)  \Bigr) = 
\end{equation}
$$
{\rm Alt}_8\Bigl(    [(3,6|1,2,5,7)_2  -
  (3,7|1,2,5,6)_2 ]
\otimes   (1,2,3,4)  \wedge (1,2,3,5)  \Bigr)
$$
using  the  five term relation
$$
 (3,6|1,2,5,7)_2 - (3,7|1,2,5,6)_2  = (3,5|1,2,6,7)_2 - (3,2|1,5,6,7)_2 + (3,1|2,5,6,7)_2 
$$
we see that the contribution of each of the last two terms is zero because 
of the skewsymmetry in $(2,3)$ and $(1,3)$.  So we get, using skew-symmetry in $(4,5)$,
$$
(\ref{ram1}) \quad = \quad {\rm Alt}_8\Bigl(     (3,4|1,2,6,7)_2  
\otimes   (1,2,3,4)  \wedge (1,2,3,5)  \Bigr)
$$

A similar consideration using the five-term relation
$$
(4,6|1,2,3,7)_2 - (4,7|1,2,3,6)_2 = (4,3|1,2,6,7)_2 -  (4,2|1,3,6,7)_2 +
  (4,1|2,3,6,7)_2 = 0 
$$
and skew-symmetry in $(1,2,3)$ gives us
$$
{\rm Alt}_8\Bigl( [2 \cdot (4,6|1,2,3,7)_2     -   3 \cdot (3,4|1,2,6,7)_2 ]  \otimes   (1,2,3,4)  \wedge (1,2,3,5)   \Bigr) = 0
$$

Using the five term relation
$$
(4,6|1,2,5,7)_2 - (4,6|1, 3,5, 6)_2 + (4,6|2,3,5,7)_2 -  (4,6|1,2,3,7)_2 +
  (4,6|1,2,3,5)_2 = 0
$$
we get
 \begin{equation} \label{ram2}
\frac{1}{3} \cdot {\rm Alt}_8\Bigl(  2 \cdot (4,6|1,2,5,7)_2       \otimes   (1,2,3,4)  \wedge (1,2,3,5)   \Bigr)  = 
 \end{equation}
$$
\frac{1}{3} \cdot {\rm Alt}_8\Bigl(   [(4,6|1,2,3,7)_2   -  (4,6|1,2,3,5)_2 ]    \otimes   (1,2,3,4)  \wedge (1,2,3,5)   \Bigr) 
$$
The last term gives zero contribution since it does not contain $7$ and $8$. 
Using the five term relation
$$
 [(4,6|1,2,3,7)_2   - [(4,7|1,2,3,6)_2  + [(4,1|2,3,6,7)_2   - [(4,2|1, 3,6, 7)_2 +[(4,3|1,2,6,7)_2 
$$
we conclude that multiplying (\ref{ram2}) by $3$ we get  
$$
{\rm Alt}_8\Bigl(    (3,4|1,2,6,7)_2     \otimes   (1,2,3,4)  \wedge (1,2,3,5)   \Bigr) 
$$

Finally,
$$
 {\rm Alt}_8\Bigl(   (4,5|1,2,6,7)_2  \otimes   (1,2,3,4)  \wedge (1,2,3,5)   \Bigr) = 0
$$
for the following reason. Write it as 
$$
\frac{1}{3}\cdot {\rm Alt}_8 \circ {\rm Alt}_{\{1,2,3\}} \Bigl( (4,5|1,2,6,7)_2  \otimes   (1,2,3,4)  \wedge (1,2,3,5)   \Bigr) 
$$
Using the five term relation for the configuration $(4,5|1,2,3,6,7)$ we can write it as 
$$
\frac{1}{3}\cdot {\rm Alt}_8 \circ {\rm Alt}_{\{6,7\}} \Bigl( (4,5|1,2,3,6)_2  \otimes   (1,2,3,4)  \wedge (1,2,3,5)   \Bigr)
$$
Each of the two terms (before taking ${\rm Alt}_8$) in this expression is zero since the first one does not contain $(7,8)$ and the second $(6,8)$. 

The lemma is proved.

\vskip 3mm \noindent
{\bf REFERENCES}
\begin{itemize}
\item[{[A]}] Aomoto K.: 
{\it Addition theorem of Abel type for hyperlogarithms},  Nagoya math. journal 88 (1982), 55-71.
\item[{[BD]}] Beilinson A.A., Deligne P.: {\it Interpr\`etation motivique de la conjecture de Zagier reliant polylogarithmes et r\'egulateurs}, Symposium in pure mathematics, 1994, vol 55, part 2, p. 97-122.
 \item[{[BMS]}] Beilinson A.A., MacPherson R., Schechtman V.V., {\it Notes on motivic cohomology}, Duke math. J., 1987, v. 54. 679-710.
\item[{[BGSV1]}] Beilinson A.A., Goncharov A.A., Schechtman V.V., Varchenko A.N.: {\it Aomoto dilogarithms, mixed Hodge structures and motivic cohomology}, the Grothendieck Festschrift, Birkhauser, vol 86, 1990, p. 135-171.
\item[{[BGSV2]}] Beilinson A.A., Goncharov A.A., Schechtman V.V., Varchenko A.N.:
{\it Projective geometry and algebraic K-theory}, Leningrad math. 
J. 2 (1991) 523-576.
\item[{[BK]}] Bloch S., Kriz I.: {\it Mixed Tate motives}, 
Annals of mathematics, 140 (1994) 557-605.
\item[{[G0]}] Goncharov A.B.: {\it The classical trilogarithm, algebraic K-theory of fields and Dedekind zeta function}, Bull AMS 29 (1991) 155-161.
\item[{[G1]}] Goncharov A.B.: {\it Geometry of configurations, polylogarithms and motivic cohomology},  Advances in Mathematics, 114, N2, (1995) 179-319.
\item[{[G2]}] Goncharov A.B.: {\it Polylogarithms and motivic Galois groups} Symposium in pure mathematics, 1994, vol 55, part 2, p.  43 - 96.
\item[{[G3]}]  Goncharov A.B.: {\it Deninger's conjecture on $L$-functions of elliptic curves at $s=3$},   J. Math. Sci., New York, 81, N 3, (1996), 2631-2656. (Special volume dedicated to Manin's 60-th birthday).
\item[{[G4]}]  Goncharov A.B.: {\it Chow polylogarithms and regulators},
Math. Res. Letters, 1995, N1, 99-114.
\item[{[G5]}]  Goncharov A.B.: {\it Multiple polylogarithms}, a paper in preparation. 
\item[{[G6]}]  Goncharov A.B.: {\it Explicit construction of
characteristic classes},  in Adv. in Soviet Mathematics,   v. 16,  N 1,   169 - 210, 1993.
\item[{[G7]}]  Goncharov A.B.: {\it Geometry of polylogarithms and regulators}, To appear.
\item[{[G8]}]  Goncharov A.B.: {\it    The generalized
Eisenstein-Kronecker series and the weight $4$ regulator}, To appear.
\item[{[G9]}]  Goncharov A.B.: {\it  Volumes of hyperbolic manifolds and mixed Tate motives}. 
Accepted for publication in Journal of AMS. 
\item[{[HM]}] Hain R., MacPherson R.: {\it Higher logarithms} 
Illinois J. math. 34 (1990) 392-475.
\item[{[HM]}] Hanamura M., MacPherson R.: {\it Geometric construction of polylogarithms}, Functional analysis on the eve of the 21st century. In honor of I.M. Gelfand, Birkhauser vol 132, (1995), 215-282.  
\end{itemize}

Dept of Mathematics, Brown University, Providence RI 02912, USA. 
 e-mail:  
sasha@math.brown.edu
\end{document}